\documentclass[a4paper,10pt]{article}
\usepackage[utf8x]{inputenc}
\usepackage{tracefnt,amsmath,tabu,array}
\usepackage{amssymb,graphicx,setspace,amsfonts,amsbsy}
\usepackage{pifont,latexsym,ifthen,amsthm,rotating,calc,textcase,booktabs,cancel,slashed}

\newtheorem{theorem}{Theorem}[section]
\newtheorem{lemma}[theorem]{Lemma}
\newtheorem{corollary}[theorem]{Corollary}
\newtheorem{proposition}[theorem]{Proposition}

\newtheorem{remark}[theorem]{Remark}
\newcommand{\filledbox}{\leavevmode
  \hbox to.77778em{%
  \hfil\vbox to.675em{\hrule width.6em height.6em}\hfil}}

\newcommand{\Rm}{{\mathbb R}}

\begin{document}
\tabulinesep=1.0mm
\title{Classification of radial non-radiative solutions to the 5D nonlinear wave equations}

\author{Liang Li, Ruipeng Shen and Chenhui Wang\\
Centre for Applied Mathematics\\
Tianjin University\\
Tianjin, China
}

\maketitle

\begin{abstract}
 In this work we classify all radial non-radiative solutions to the 5D nonlinear wave equations with a wide range of energy critical nonlinearity. We show that such a solution always comes with two characteristic numbers. These characteristic numbers can be determined by either the radiation profile of the initial data or the asymptotic behaviour of the solution. In addition, two radial weakly non-radiative solutions with the same characteristic numbers must coincide with each other in the overlap part of their exterior regions. Finally we give a few applications of our theory on the global behaviours of solutions to the nonlinear wave equations.
\end{abstract}

\section{Introduction}

\paragraph{Assumptions} We consider the Cauchy problem of 5D wave equation 
\[
 \left\{\begin{array}{ll} \partial_t^2 u - \Delta u = F(t,x,u);\\ (u,u_t)|_{t=0} = (u_0,u_1).  \end{array}\right. \quad (CP1)
\]
In this article the nonlinear term $F(t,x,u)$ always satisfies the following assumptions 
\begin{itemize} 
\item The function $F(t,x,u)$ is a radial function of $x$;
\item The inequalities below hold for a nonnegative constant $\gamma$. 
\begin{align*}
&|F(t,x,u)|  \leq \gamma |u|^{7/3}; &  |F(t,x,u_1)-F(t,x,u_2)| \leq \gamma (|u_1|^{4/3}+|u_2|^{4/3})|u_1-u_2|;&
\end{align*}
\end{itemize}
Some of our results also assume the symmetric properties 
\begin{align*}
 & F(-t,x,u) = F(t,x,u); & & F(t,x,-u)  = -F(t,x,u). & &(AS)&
\end{align*}
Examples of this kind of equations include the focusing/defocusing equations with a power-type nonlinearity 
\[
 \partial_t^2 u - \Delta u = \pm |u|^{4/3} u,
\]
which are extensively studied in the past decades. Please see \cite{kenig}, for example. 

\paragraph{Topic of this work} The channel of energy method plays an important role in the study of nonlinear wave equation in recent years. The application of this method includes the proof of soliton resolution conjecture of energy critical wave equation with radial data in all odd dimensions $d\geq 3$ by Duyckaerts-Kenig-Merle \cite{se, oddhigh} and the conditional scattering of solutions in the energy super or sub-critical case by Duyckaerts-Kenig-Merle \cite{dkm2} and Shen \cite{shen2}. The study of non-radiative solutions is a major topic in the theory of channel of energy method. In this work we classify all radial non-radiative solutions to the 5D energy critical wave equation. We also give asymptotic behaviour of these solution in details. Before we give our main result, we first introduce some necessary backgrounds and a few basic conceptions.  

\paragraph{Local theory in the whole space} The equation (CP1) is well-posed in the space $\dot{H}^1 \times L^2(\Rm^5)$. The idea is to combine suitable Strichartz estimates (see \cite{strichartz}) and a standard fixed-point argument. One may follow a similar argument to those given in Kapitanski \cite{loc1}, Kenig-Merle \cite{kenig} and Lindblad-Sogge \cite{ls}. Given an time interval $I$ containing zero, we call $u$ a solution to (CP1) in the time interval $I$ if and only if 
\begin{itemize}
 \item Given any finite bounded interval $J \subset I$, we have\footnote{Some literatures use other space-time norms, but all these definitions are actually equivalent by the Strichartz estimates.} $\|u\|_{Y(J)} < +\infty$;
 \item The solution satisfies 
 \[
  u = \mathbf{S}_L (u_0,u_1) + \int_0^t \frac{\sin (t-t')\sqrt{-\Delta}}{\sqrt{-\Delta}} F(t', x, u) dt', \qquad t\in I.
 \]
\end{itemize}
Here $\mathbf{S}_L (u_0,u_1)$ is the solution to the homogenous linear wave equation $\partial_t^2 u - \Delta u = 0$ with initial data $(u_0,u_1)$ and the $Y$ space-time norm is defined by 
\[
 \|u\|_{Y(I)}  = \|u\|_{L_t^{7/3} L_x^{14/3} (I \times \Rm^5)} = \left(\int_I \left(\int_{\Rm^5} |u(x,t)|^{14/3}\right)^{1/2}dt\right)^{3/7}.
\]
We also need to use the following norm in this work
\[
  \|F\|_{Z(I)}  = \|F\|_{L_t^{1} L_x^{2} (I \times \Rm^5)} = \int_I \left(\int_{\Rm^5} |F(x,t)|^2\right)^{1/2}dt.
\]

\paragraph{Exterior solutions} Let $I$ be a time interval containing $0$. We say $u$ defined in the exterior region $\{(x,t):t\in I, |x|>R+|t|\}$ is an exterior solution to (CP1) with initial data $(u_0,u_1)\in \dot{H}^1 \times L^2(\Rm^5)$ , if and only if 
\begin{itemize}
 \item Given any finite bounded interval $J \subset I$, we have $\|\chi_R u\|_{Y(J)} < +\infty$;
 \item The solution satisfies 
 \[
  u = \mathbf{S}_L (u_0,u_1) + \int_0^t \frac{\sin (t-t')\sqrt{-\Delta}}{\sqrt{-\Delta}} [\chi_R F(t', x, u)] dt', \qquad t\in I, |x|>R+|t|
 \]
\end{itemize}
Here the factor $\chi_R$ is the characteristic function of the region $\Omega_R = \{(x,t)\in \Rm^5 \times \Rm: |x|>|t|+R\}$. The notations $\Omega_R$ and $\chi_R$ will be used throughout this work. We put the factor $\chi_R$ in $\chi_R u$ and $\chi_R F(t,x,u)$ to guarantee that these functions are defined globally in $\Rm^5$ so we may apply the $Y$ norm and the operator defined by Fourier multipliers conveniently. We recall the Strichartz estimates (details are given in Section \ref{sec: preliminary})
\[
 \|u\|_{C(I; \dot{H}^1 \times L^2)} + \|u\|_{Y(I)} \lesssim \|(u_0,u_1)\|_{\dot{H}^1 \times L^2} + \|(\partial_t^2 - \Delta)u\|_{Z(I)}
\]
 and observe the inequalities 
\begin{align*}
 \|\chi_R  F(t, x, u)\|_{Y(I)} &\lesssim \|\chi_R u\|_{Y(I)}^{7/3}; \\
 \|\chi_R F(t,x,u_1) - \chi_R F(t,x,u_2)\|_{Y(I)} & \lesssim \left(\|\chi_R u_1\|_{Y(I)}^{4/3} + \|\chi_R u_2\|_{Y(I)}^{4/3}\right) \|\chi_R (u_1-u_2)\|_{Y(I)}.
\end{align*}
A standard fixed-point argument then gives the well-posedness of this Cauchy problem. Please note that the exterior solution depends on the values of the initial data in the exterior region $\{x: |x|>R\}$ only.  

\paragraph{Restriction of solutions} Let $R\geq 0$. If $u$ is a solution to (CP1) defined in $\Rm^5 \times I$, or an exterior solution defined in $\{(x,t): t\in I, |x|>r+|t|\}$ with $0 \leq r < R$. Then the restriction of $u$ in the exterior region $\{(x,t): t\in I, |x|>R+|t|\}$ is an exterior solution to (CP1). We call it the $R$-restriction of $u$. 

\paragraph{Global extension} If $u$ is an exterior solution defined in the exterior region $\{(x,t):t\in I, |x|>R+|t|\}$, then we may consider the solution $\tilde{u}$ to the linear wave equation 
\[
 \left\{\begin{array}{ll} \partial_t^2 \tilde{u} - \Delta \tilde{u} = \chi_R F(t,x,u), & (x,t) \in \Rm^5 \times I;\\ (\tilde{u}, \tilde{u}_t)|_{t=0} = (u_0,u_1). & \end{array} \right.
\]
By the Strichartz estimates, we have $(\tilde{u}, \tilde{u}_t) \in C(I; \dot{H}^1 \times L^2(\Rm^5))$. Finite speed of propagation shows that $\tilde{u}$ coincides with $u$ in the exterior region $\{(x,t): t\in I, |x|>|t|+R\}$. Thus $\tilde{u}$ also solves the non-linear wave equation 
\begin{equation}
 \partial_t^2 \tilde{u} - \Delta \tilde{u} = \chi_R F(t,x,\tilde{u}) \label{localized equation}
\end{equation}
We call $\tilde{u}$ defined above an global extension of $u$ in $\Rm^5 \times I$. Conversely, if $\tilde{u}$ is a solution to \eqref{localized equation} defined in $\Rm^5 \times I$, then its restriction in the exterior region $\{(x,t): t\in I, |x|>|t|+R\}$ is an exterior solution to (CP1). Please note that unlike the exterior solutions, the extensions of exterior solutions depend on the values of initial data $(u_0,u_1)$ in the interior region $\{x: |x|<R\}$. If the values of initial data are not specified for an exterior solution, we have to define them in some way before we may consider the global extensions. 

\paragraph{Non-radiative solutions} Let $R\geq 0$. We call an exterior solution $u$ defined in the exterior region $\Omega_R = \{(x,t): |x|>|t|+R\}$ an $R$-weakly non-radiative solution if and only if 
\[
 \lim_{t\rightarrow \pm \infty} \int_{|x|>R+|t|} |\nabla_{t,x} u (x,t)|^2 dx = 0. 
\]
If $u$ is defined in a larger region, for example, the exterior region $\Omega_r$ with $0 \leq r < R$, or the whole space $\Rm^5 \times \Rm$, then we also call it an $R$-weakly non-radiative solution if it satisfies the limit above. In particular, we call the solutions non-radiative if the limit holds for $R=0$. If it is unnecessary to mention the radius $R$, we may also call these solutions weakly non-radiative solutions. 

\paragraph{Radiation fields} Generally speaking, radiation fields describe the asymptotic behaviour of linear free waves. The history of radiation fields are more than 50 years long. Please see, Friedlander \cite{radiation1, radiation2} for example. The following version of radiation fields comes from Duyckaerts-Kenig-Merle \cite{dkm3} and plays an important role in the discussion of asymptotic behaviours of non-linear wave equations. The radiation field is the most important tool in this work. 
\begin{theorem}[Radiation field] \label{radiation}
Assume that $d\geq 3$ and let $u$ be a solution to the free wave equation $\partial_t^2 u - \Delta u = 0$ with initial data $(u_0,u_1) \in \dot{H}^1 \times L^2(\Rm^d)$. Then ($u_r$ is the derivative in the radial direction)
\[
 \lim_{t\rightarrow \pm \infty} \int_{\Rm^d} \left(|\nabla u(x,t)|^2 - |u_r(x,t)|^2 + \frac{|u(x,t)|^2}{|x|^2}\right) dx = 0
\]
 and there exist two functions $G_\pm \in L^2(\Rm \times \mathbb{S}^{d-1})$ so that
\begin{align*}
 \lim_{t\rightarrow \pm\infty} \int_0^\infty \int_{\mathbb{S}^{d-1}} \left|r^{\frac{d-1}{2}} \partial_t u(r\theta, t) - G_\pm (r\mp t, \theta)\right|^2 d\theta dr &= 0;\\
 \lim_{t\rightarrow \pm\infty} \int_0^\infty \int_{\mathbb{S}^{d-1}} \left|r^{\frac{d-1}{2}} \partial_r u(r\theta, t) \pm G_\pm (r\mp t, \theta)\right|^2 d\theta dr & = 0.
\end{align*}
In addition, the maps $(u_0,u_1) \rightarrow \sqrt{2} G_\pm$ are bijective isometries from $\dot{H}^1 \times L^2(\Rm^d)$ to $L^2 (\Rm \times \mathbb{S}^{d-1})$. 
\end{theorem}
\noindent We call $G_\pm$ the radiation profile of $u$, or equivalently, its initial data $(u_0,u_1)$. In the 5-dimensional case, the radiation profiles $G^\pm$ satisfies $G_+(s,\theta) = G_-(-s,-\theta)$, as shown in C\^{o}te-Laurent \cite{newradiation} and Li-Shen-Wei \cite{shenradiation}. Therefore a linear free wave is $R$-weakly non-radiative if and only if its radiation profiles are compactly supported in the region $[-R,R]\times \mathbb{S}^4$. The one-to-one map between radiation profiles and linear free waves (or equivalently, their initial data) can be explicitly given in many different ways. The following formula in the 5D case comes from Li-Shen-Wei \cite{shenradiation}. 
\[
 u(x,t) = \frac{1}{4\pi^2} \int_{\mathbb{S}^4} G'_-(x\cdot \omega +t , \omega) d\omega. 
\]
An explicit formula in term of Fourier transforms can also be found in a recent work C\^{o}te-Laurent \cite{newradiation}. In this work, if we mention the radiation profile of a free wave, or equivalently its initial data, then we mean the radiation profile in the negative time direction unless specified otherwise. It is clear that the free wave is a radial function of $x$ if and only if the radiation profile $G$ is independent of $\theta$. In this work we only consider radial solutions, thus the radiation profiles are always viewed as one-variable functions $G(s)$. 

\paragraph{Goal of this work} We discuss the classification of all radial weakly non-radiative solutions to (CP1) in this work. Duyckaerts-Kenig-Merle \cite{oddtool} describes the asymptotic behaviours of all radial weakly non-radiative solutions to the focusing equation $\partial_t^2 u - \Delta u = +|u|^{4/(d-2)}u$ in all odd dimensions $d\geq 5$. For simplicity we only give the statement in 5-dimensional case here, although higher dimensions are similar. 
\begin{theorem} \label{Kenig asymptotic} 
 Let $\Theta_1 = (t|x|^{-3}, |x|^{-3})$ and $\Theta_2 = (|x|^{-3}, 0)$. If $u$ is a radial weakly non-radiative solution to $\partial_t^2 u - \Delta u = + |u|^{4/3} u$, then there exists $R_1 \gg 1$, $k \in \{1,2\}$, $l \in \Rm$ with $l \neq 0$ if $k =1$, such that for all $t \in \Rm$, 
 \[ 
  \|(u(\cdot,t), u_t(\cdot,t))- l \Theta_k (\cdot,t) \|_{\dot{H}^1 \times L^2(\{x:|x|>R\})}\lesssim \max\{R^{-\frac{7}{3}(k-\frac{1}{2})}, R^{-(k+\frac{1}{2})}\}, \quad \forall R > |t|+R_1.  
 \]
 Furthermore, if $k=2$, then the solution $u$ can be given explicitly 
 \[
  u(x,t) = \left(\frac{1}{\lambda} + \frac{\lambda |x|^2}{15}\right)^{-3/2}, \qquad |x|>R_1+|t|.
 \]
\end{theorem}
\noindent Generally speaking, a weakly non-radiative solution share the same asymptotic behaviour as $l \Theta_m (x,t)$. This is quite reasonable since it has been proved in Kenig et al \cite{channel} that all radial $R$-weakly non-radiative solutions to the homogenous linear wave equation are exactly those in the linear space 
\[
  \hbox{Span} \{t|x|^{-3}, |x|^{-3}\}.
\]
In this work we are trying to answer the following questions: Given $l \neq 0$, does there exists a radial weakly non-radiative solution whose asymptotic behaviour is the same as that of $lt|x|^{-3}$? If there does, how many solutions of this kind are there? Can we further classify these solutions? Next we give our main result.

\begin{theorem} \label{main 1} 
 We consider radial weakly non-radiative solutions to (CP1). 
 \begin{itemize}
  \item [(a)] If $u$ is such a solution, then the radiation profile $G$ of its initial data satisfies 
  \begin{align*}
   & \|G\|_{L^2\{s:|s|>r\}} \lesssim r^{-7/6}, \quad r\gg 1;& &\Rightarrow& &G \in L^1(\Rm).&
  \end{align*}
  We call the number
  \begin{align*}
    \alpha = -\int_{-\infty}^\infty G(s) ds
  \end{align*}
  the first characteristic number of $u$. The first characteristic number can also be characterized by the asymptotic behaviour of the initial data 
  \begin{equation} \label{asymptotic first order}
   \|(u_0,u_1) - (0,\alpha |x|^{-3})\|_{\dot{H}^1 \times L^2(\{x: |x|>r\})} \lesssim r^{-7/6}, \qquad r\gg 1.\\
  \end{equation} 
  In addition, given any $\alpha$, there exists at least one radial weakly non-radiative solution to (CP1), whose first characteristic number is exactly $\alpha$. Given $\alpha \in \Rm$, we choose such a solution and call it $u^\alpha$. We also use the notation $G_\alpha$ for the radiation profile of its initial data.
  \item [(b)] Given $\alpha \in \Rm$ and $u^\alpha$, if $u$ is a radial weakly non-radiative solution to (CP1) whose first characteristic number is $\alpha$, then the radiation profile $G$ of its initial data $(u_0,u_1)$ satisfies 
  \begin{align*} 
   &\|G-G_\alpha\|_{L^2(\{s: |s|>r\})} \lesssim r^{-13/6}, \quad r\gg 1;& &\Rightarrow& &s[G(s)-G_\alpha (s)] \in L^1 (\Rm).&
  \end{align*}
  We call the number 
  \[
   \int_{-\infty}^{+\infty} s[G(s)-G_\alpha(s)] ds
  \]
  the second characteristic number (with respect to $u^\alpha$). The second characteristic number can also be characterized by the asymptotic behaviour of $u$. 
  \begin{equation} \label{asymptotic second order}
   \left\|(u_0,u_1)-(u^\alpha(\cdot,0), u_t^\alpha(\cdot,0)) - (\beta |x|^{-3},0)\right\|_{\dot{H}^1 \times L^2(\{x:|x|>r\})} \lesssim r^{-13/6}, \qquad r\gg 1.
  \end{equation}
  In addition, given any $\alpha, \beta \in \Rm$, there exists a radial weakly non-radiative solution to (CP1) whose first and second characteristic numbers are exactly $\alpha$ and $\beta$. 
  \item[(c)] Let $u$ and $v$ be two radial $R$-weakly non-radiative solutions to (CP1) so that they share the same first and second characteristic numbers. Then they must be identical to each other.
  \[
   u(x,t) = v(x,t), \qquad \hbox{if}\; |x|>|t|+R. 
  \]
 \end{itemize}
\end{theorem}

\begin{remark}
Let $u$ be a radial $R$-weakly non-radiative solution to (CP1). Its characteristic numbers do not depend on the value of initial data $(u_0,u_1)$ in the interior region $\{x: |x|<R\}$, although the radiation profile $G$ does depend on these values. It is natural because the exterior solution does not depend on the value of initial data in the interior region. This independence can be verified in either of the following two ways.
\begin{itemize}
 \item By the definition of radiation profile given in Theorem \ref{radiation}, the values of $G(s)$ with $|s|>r$ only depend on the value of linear free wave in the exterior region $\Omega_R$, thus are completely determined by the value of $(u_0,u_1)$ in $\{x: |x|>R\}$, by finite speed of propagation. In addition, the integrals
\begin{align*}
 &\int_{-R}^R G(s) ds,& &\int_{-R}^R sG(s) ds&
\end{align*}
are also determined by the value of $(v_0,v_1)$ in the region $\{x: |x|>R\}$, by Lemma \ref{u by G}. 
\item The characteristic numbers $\alpha, \beta$ can be uniquely determined by the asymptotic behaviour of initial data. In fact, the pair $(\alpha,\beta)$ satisfying \eqref{asymptotic first order} and \eqref{asymptotic second order} is unique since we have 
\begin{align*} 
  &\||x|^{-3}\|_{L^2(\{x: |x|>r\})} \simeq r^{-1/2};& &\||x|^{-3}\|_{\dot{H}^1(\{x: |x|>r\})} \simeq r^{-3/2}.&
 \end{align*}
\end{itemize}
\end{remark}

\begin{remark} 
Given $\alpha \in \Rm\setminus \{0\}$, we usually choose $u^\alpha$ to be the weakly non-radiative solution constructed in Lemma \ref{existence of solution with first alpha}. In this work all the second characteristic numbers are chosen with respect to these solutions $u^\alpha$ unless specified otherwise. If the nonlinear term $F(t,x,u)$ satisfies the additional symmetric assumption (AS), then the solutions $u^\alpha$ mentioned above are natural choices because they also satisfies similar symmetric property $u^\alpha (x,-t) = - u^\alpha(x,t)$. In this case, the second characteristic number can be determined conveniently by the integral of radiation profile $G$ or the asymptotic behaviour of initial data, without the explicit appearance of $u^\alpha$ or $G_\alpha$. More precisely we have
\[
 \beta = \hbox{p.v.} \int_{-\infty}^{+\infty} s G(s) ds \doteq \lim_{r\rightarrow +\infty} \int_{-r}^r sG(s) ds;
\]
and 
\[
  \|u_0 - \beta |x|^{-3}\|_{\dot{H}^1(\{x: |x|>r\})} \lesssim r^{-13/6}, \qquad r\gg 1.
\]
\end{remark}

\begin{remark}
Our main theorem implies that a radial weakly non-radiative solution is essentially determined by its two characteristic numbers. Since the characteristic numbers can be determined by the asymptotic behaviour of solutions, a radial weakly non-radiative solution is also essentially determined by its values in an exterior region $\Omega_r$ with an arbitrarily large radius $r$.  More precisely, if $u$ and $v$ are both radial weakly non-non-radiative solutions to (CP1), defined in the exterior regions $\Omega_{R_1}$ and $\Omega_{R_2}$ respectively, then the following three statements are equivalent to each other:
\begin{itemize}
 \item $u(x,t) = v(x,t)$ holds in the overlap part of their exterior regions $\Omega_{\max\{R_1,R_2\}} = \Omega_{R_1} \cap \Omega_{R_2}$;
 \item $u(x,t) = v(x,t)$ holds in any exterior region $\Omega_r$ with $r \geq \max\{R_1,R_2\}$;
 \item $u$ and $v$ share the same first and second characteristic numbers. 
\end{itemize}
\end{remark}

\begin{remark}
The number $l$ in Duyckaerts-Kenig-Merle's Theorem \ref{Kenig asymptotic} is actually the first nonzero characteristic number, unless the solution is identical to zero.   
\end{remark}

\begin{remark}
 A recent paper Collot-Duyckaerts-Kenig-Merle \cite{kenignparameters} classifies all radial non-radiative solutions to the energy critical wave equations $\partial_t^2 u - \Delta u = |u|^{4/(d-2)} u$ for all odd dimensions $d \geq 3$ in a different method. They show that all the radial non-radiative solutions form a family with $(d-1)/2$ parameters. But their choice of parameters are different from ours. They first choose polynomials $p_k$'s so that all the radial $R$-weakly non-radiative solutions to the free wave equation are exactly the linear combinations ($c_k$ are constants)
 \[
   \sum_{k=0}^{(d-3)/2} \frac{c_k}{|x|^{d-2-k}} p_k \left(\frac{t}{|x|}\right),
 \]
and let $\vec{a}(\{c_k\})$ be the corresponding initial data. For convenience the notation $\mathcal{P}$ is used for the linear space of all such initial data. 
If $u$ be a radial $R$-weakly non-radiative solution to the non-linear equation so that the energy of initial data in the exterior region $\{x: |x|>R\}$ is sufficiently small, then the parameters are chose to be the coefficients $c_k$'s so that 
\[
 \vec{a}(\{c_k\}) = \mathbf{\Pi}_{\mathcal{P}} (u(\cdot,0), u_t(\cdot,0)).
\]
Here $\mathbf{\Pi}_{\mathcal{P}}$ is the orthogonal projection of the space $\dot{H}^1\times L^2(\{x: |x|>R\})$ onto its finite dimensional space $\mathcal{P}$. Please note that the values of their parameters may depend on the radius $R$. Our characteristic numbers, however, do not depend on the choice of radius. 
\end{remark}

\paragraph{Structure of this work} We first give a few preliminary results and technical lemmata in Section 2. The existence and properties of first and second characteristic numbers are proved in Section 3 and 4, respectively. The final section is devoted to further analysis and applications of characteristic numbers.

\paragraph{Notations} In this work the notation $A \lesssim B$ means that there exists a constant $c>0$, so that the inequality $A \leq cB$ always holds.  We may also add subscript(s) to the symbol $\lesssim$ to emphasize that the constant $c$ depends on the subscript(s) but nothing else. In particular, $\lesssim_1$ means that the constant $c$ is an absolute constant. 

\section{Preliminary Results} \label{sec: preliminary}

\paragraph{Strichartz estimates} We first recall the generalized Strichartz estimates. Please see Proposition 3.1 in Ginibre-Velo \cite{strichartz}. Here we use the Sobolev version in dimension $5$. 

\begin{proposition}[Strichartz estimates]\label{Strichartz estimates} 
 Let $2\leq q_1,q_2 \leq \infty$, $2\leq r_1,r_2 < \infty$ and $\rho_1,\rho_2,s\in \Rm$ be constants with
 \begin{align*}
  &\frac{1}{q_i} + \frac{2}{r_i} \leq 1,& &(q_i,r_i)\neq \left(2,4\right),& &i=1,2;& \\
  &\frac{1}{q_1} + \frac{5}{r_1} = \frac{5}{2} + \rho_1 - s;& &\frac{1}{q_2} + \frac{5}{r_2} = \frac{3}{2} + \rho_2 +s.&
 \end{align*}
 Assume that $u$ is the solution to the linear wave equation
\[
 \left\{\begin{array}{ll} \partial_t u - \Delta u = F(x,t), & (x,t) \in \Rm^5 \times [0,T];\\
 u|_{t=0} = u_0 \in \dot{H}^s; & \\
 \partial_t u|_{t=0} = u_1 \in \dot{H}^{s-1}. &
 \end{array}\right.
\]
Then we have
\begin{align*}
 \left\|\left(u(\cdot,T), \partial_t u(\cdot,T)\right)\right\|_{\dot{H}^s \times \dot{H}^{s-1}} & +\|D_x^{\rho_1} u\|_{L^{q_1} L^{r_1}([0,T]\times \Rm^5)} \\
 & \leq C\left(\left\|(u_0,u_1)\right\|_{\dot{H}^s \times \dot{H}^{s-1}} + \left\|D_x^{-\rho_2} F(x,t) \right\|_{L^{\bar{q}_2} L^{\bar{r}_2} ([0,T]\times \Rm^5)}\right).
\end{align*}
Here $\bar{q}_2$ and $\bar{r}_2$ satisfy $1/q_2 + 1/\bar{q}_2 = 1$, $1/r_2 + 1/\bar{r}_2 = 1$. The constant $C$ does not depend on $T$ or $u$. 
\end{proposition}
\noindent In particular, we may choose $(q_1,r_1,\rho_1) = (7/3,14/3,0)$ and $(q_2,r_2,\rho_2) = (\infty, 2, 0)$. The solution to the linear wave equation above satisfies 
\begin{align*}
 \left\|\left(u(\cdot,T), \partial_t u(\cdot,T)\right)\right\|_{\dot{H}^s \times \dot{H}^{s-1}} & +\|u\|_{L^{7/3} L^{14/3}([0,T]\times \Rm^5)} \\
 & \leq C\left(\left\|(u_0,u_1)\right\|_{\dot{H}^s \times \dot{H}^{s-1}} + \left\|F(x,t) \right\|_{L^{1} L^{2} ([0,T]\times \Rm^5)}\right).
\end{align*}

\begin{lemma} \label{u by G}
 Assume that $u$ is a radial free wave with initial data $(u_0,u_1) \in \dot{H}^1 \times L^2(\Rm^5)$. Let $G(s) \in L^2(\Rm)$ be the radiation profile in the negative time direction associated to $u$. Then we have 
 \[
  u(r,t) = \frac{1}{r^3}\int_{t-r}^{t+r} (s-t) G(s) ds. 
 \]
 In particular, the initial data are given by
 \begin{align*}
  &u_0(r,t) =  \frac{1}{r^3}\int_{-r}^{+r} s G(s) ds; & &u_1(r,t) = \frac{G(r)+G(-r)}{r^2} -\frac{1}{r^3} \int_{-r}^r G(s) ds. &
 \end{align*}
\end{lemma} 
\begin{proof}
 Without loss of generality, we assume that $G$ is smooth and compactly supported. The general case follows standard smooth approximation techniques. We first recall the following formula of free wave $u$ in term of the radiation profile $G(s,\omega)$ given in Remark 2.3 of \cite{shenradiation}. 
\[
 u(x,t) = \frac{1}{4\pi^2} \int_{\mathbb{S}^4} G'(x\cdot \omega +t , \omega) d\omega. 
\]
Here $G'$ is the derivative of $G$ with respect to the first variable $s$. We next utilize the radial assumption and integrate ($\omega = (\omega_1, \omega_2, \cdots, \omega_5) \in \mathbb{S}^4$)
\[
 u(r,t) = u(r,0,0,0,0,t) = \frac{1}{4\pi^2} \int_{\mathbb{S}^4} G'(r\omega_1+t) d\omega = \frac{1}{2} \int_{-1}^1 G'(r\omega_1+t) (1-\omega_1^2) d\omega_1. 
\]
We then integrate by parts and obtain 
\[
 u(r,t) = \frac{1}{r} \int_{-1}^1 G(r\omega_1+t) \omega_1 d\omega_1 = \frac{1}{r^3} \int_{t-r}^{t+r} (s-t) G(s) ds. 
\]
A straight forward calculation then gives the formula of $(u_0,u_1)$. 
\end{proof}

\begin{corollary} \label{outward decay}
 Let $u$ be a radial $R$-weakly non-radiative solution to the free wave equation $\partial_t^2 u - \Delta u = 0$ with radiation profile $G$. Then we have 
\[
 \|(u_0,u_1)\|_{\dot{H}^1 \times L^2(\{x: |x|>r\})} + \|\chi_r u\|_{Y(\Rm)} \lesssim_1 r^{-3/2} \left|\int_{-R}^R sG(s) ds\right| + r^{-1/2} \left|\int_{-R}^R G(s) ds\right|, \qquad r\geq R.
\]
In particular we have 
\[
 \|\chi_r u\|_{Y(\Rm)} \lesssim_1 (R/r)^{1/2} \|G\|_{L^2(\Rm)}, \qquad r\geq R. 
\]
\end{corollary}
\begin{proof}
We conduct a straight forward calculation. Since $u$ is $R$-weakly non-radiative free wave, we have $G(s) = 0$ for $s>R$. Therefore we may apply the formula given in Lemma \ref{u by G} to obtain that if $r > R$, then
\begin{align*}
 &u_0(r,t) =  \frac{1}{r^3}\int_{-R}^{+R} s G(s) ds; & &u_1(r,t) = -\frac{1}{r^3} \int_{-R}^R G(s) ds. &
\end{align*}
A straight forward calculation then gives the upper bound of $\|(u_0,u_1)\|_{\dot{H}^1 \times L^2(\{x: |x|>r\})}$. The upper bound of $\|\chi_r u\|_{Y(\Rm)}$ then follows Strichartz estimates and finite speed of propagation. The last inequality immediately follows the Cauchy-Schwartz 
\begin{align*}
 &\left|\int_{-R}^R sG(s) ds\right| \lesssim_1 R^{3/2} \|G\|_{L^2(\Rm)};& &\left|\int_{-R}^R G(s) ds\right| \lesssim_1 R^{1/2} \|G\|_{L^2(\Rm)}.&
\end{align*}
\end{proof}

\begin{corollary} \label{outward decay 2}
 Let $(u_0,u_1)\in \dot{H}^1 \times L^2(\Rm^5)$ be radial initial data whose radiation profile is $G$. Then given $r>0$, we have 
 \[
  \|(u_0,u_1)\|_{\dot{H}^1 \times L^2(\{x: |x|>r\})} \lesssim_1 r^{-3/2} \left|\int_{-r}^r sG(s) ds\right| + r^{-1/2} \left|\int_{-r}^r G(s) ds\right| + \|G\|_{L^2(\{s: |s|>r\})}. 
 \]
 In addition, the linear free wave $u = \mathbf{S}_L(u_0,u_1)$ satisfies
 \[
  \|\chi_r u\|_{Y(\Rm)} \lesssim_1 r^{-3/2} \left|\int_{-r}^r sG(s) ds\right| + r^{-1/2} \left|\int_{-r}^r G(s) ds\right| + \|G\|_{L^2(\{s: |s|>r\})}. 
 \]
\end{corollary}
\begin{proof}
We first split $G$ into two parts
  \begin{align*}
   &G_1(s) = \left\{\begin{array}{ll} G(s), & |s|\leq r; \\ 0, & |s| > r; \end{array}\right. &  &G_2(s) = \left\{\begin{array}{ll} 0, & |s|\leq r; \\ G(s), & |s| > r; \end{array}\right. & 
  \end{align*}
 and write $(u_0,u_1) = (u_0^1, u_1^1) + (u_0^2, u_1^2)$ accordingly. We have 
 \[
  \|(u_0,u_1)\|_{\dot{H}^1 \times L^2(\{x: |x|>r\})} \leq \|(u_0^1,u_1^1)\|_{\dot{H}^1 \times L^2(\{x: |x|>r\})} + \|(u_0^2,u_1^2)\|_{\dot{H}^1 \times L^2(\{x: |x|>r\})}.
 \]
 We apply Corollary \ref{outward decay} on $(u_0^1, u_1^1)$ and the isometric identity on $(u_0^2,u_1^2)$ 
 \begin{align*}
  \|(u_0^1,u_1^1)\|_{\dot{H}^1 \times L^2(\{x: |x|>r\})} & \lesssim_1 r^{-3/2} \left|\int_{-r}^r sG(s) ds\right| + r^{-1/2} \left|\int_{-r}^r G(s) ds\right|;\\
  \|(u_0^2,u_1^2)\|_{\dot{H}^1 \times L^2(\{x: |x|>r\})} & \lesssim_1 \|G_2\|_{L^2(\Rm)} = \|G\|_{L^2(\{s: |s|>r\})}.
 \end{align*}
 Plugging these upper bounds in the inequality above, we obtain the upper bound of $L^2$ norm. The upper bound of $\|\chi_r u\|_{Y(\Rm)}$ follows the Strichartz estimates and finite speed of propagation. 
\end{proof}
\begin{lemma} \label{Y linear to nonlinear}
 Let $(u_0,u_1)\in \dot{H}^1 \times L^2(\Rm^5)$ be initial data so that $\|\chi_r \mathbf{S}_L (u_0,u_1)\|_{Y(\Rm)}$ is sufficiently small, then there exists a unique exterior solution $u$ to (CP1) defined in $\Omega_r$ with initial data $(u_0,u_1)$. We also have
 \[
  \|\chi_r u\|_{Y(\Rm)} \leq 2 \|\chi_r \mathbf{S}_L (u_0,u_1)\|_{Y(\Rm)}. 
 \] 
 In addition, if $(v_0,v_1)$ are initial data satisfying the same assumption, then the corresponding exterior solution $v$ satisfies 
 \[
  \|\chi_r (u-v)\|_{Y(I)} \leq 2  \|\chi_r \mathbf{S}_L (u_0-v_0,u_1-v_1)\|_{Y(\Rm)}. 
 \]
\end{lemma}
\begin{proof}
The existence, uniqueness and upper bound of $u$ immediately follows a standard fixed-point argument. Given two such pairs of initial data, we may apply the Strichartz estimate and obtain 
\begin{align*}
 \|\chi_r (u-v)\|_{Y(\Rm)} &\leq \|\chi_r \mathbf{S}_L(u_0-v_0,u_1-v_1)\|_{Y(\Rm)} + C_1 \|\chi_r F(t,x,u) - \chi_r F(t,x,v)\|_{Z(\Rm)}\\
 & \leq \|\chi_r \mathbf{S}_L(u_0-v_0,u_1-v_1)\|_{Y(\Rm)} + C_2 \left(\|\chi_r u\|_{Y(\Rm)}^{4/3} + \|\chi_r v\|_{Y(\Rm)}^{4/3} \right)\|\chi_r (u-v)\|_{Y(\Rm)}.
\end{align*}
Here $C_1, C_2$ are constants determined solely by $\gamma$. If $\|\chi_r \mathbf{S}_L (u_0,u_1)\|_{Y(\Rm)}$ and $\|\chi_r \mathbf{S}_L (v_0,v_1)\|_{Y(\Rm)}$ are sufficiently small, we have 
\[
 C_2 \left(\|\chi_r u\|_{Y(\Rm)}^{4/3} + \|\chi_r v\|_{Y(\Rm)}^{4/3} \right) < 1/2,
\]
thus 
\[
 \|\chi_r (u-v)\|_{Y(\Rm)}  \leq 2\|\chi_r \mathbf{S}_L (u_0-v_0,u_1-v_1)\|_{Y(\Rm)}. 
\]
\end{proof}

\begin{lemma} [Radiation fields of inhomogeneous equation]\label{scatter profile of nonlinear part}
 Let $u$ be a radial solution to the linear wave equation
 \[
  \left\{\begin{array}{ll} \partial_t^2 u - \Delta u = F(t,x); & (x,t)\in \Rm^5 \times \Rm; \\
  (u,u_t)|_{t=0} = (0,0). & \end{array} \right.
 \]
 If $F\in Z(\Rm)$ is a radial function, then there exists $G^\pm \in L^2(\Rm)$ so that 
 \begin{align*}
  \lim_{t\rightarrow +\infty} \int_0^\infty \left(\left|G^+(r-t) - r^2 u_t (r, t)\right|^2 + \left|G^+(r-t) + r^2 u_r (r, t)\right|^2\right) dr & = 0; \\
  \lim_{t\rightarrow -\infty} \int_0^\infty \left(\left|G^-(r+t) - r^2 u_t(r,t)\right|^2 +  \left|G^-(r+t) - r^2 u_r(r,t)\right|^2\right)dr & = 0.
 \end{align*}
 In addition, we have the upper bounds ($R \geq 0$)
 \begin{align*}
  &\|G^-\|_{L^2([R,+\infty))} \lesssim_1 \|\chi_R F\|_{Z((-\infty,0])};&  &\|G^+\|_{L^2([R,+\infty))} \lesssim_1 \|\chi_R F\|_{Z([0,+\infty))}.&
 \end{align*}
\end{lemma}
\begin{proof}
We may apply the Strichartz estimates and obtain 
\begin{align*}
&  \lim_{t_1, t_2 \rightarrow +\infty} \left\|\mathbf{S}_L(-t_1) \begin{pmatrix} u(\cdot,t_1)\\ u_t(\cdot,t_1)\end{pmatrix} - \mathbf{S}_L(-t_2) \begin{pmatrix} u(\cdot,t_2)\\ u_t(\cdot,t_2)\end{pmatrix}\right\|_{\dot{H}^1 \times L^2(\Rm^5)} \\
= & \lim_{t_1, t_2 \rightarrow +\infty} \left\|\mathbf{S}_L(t_2-t_1) \begin{pmatrix} u(\cdot,t_1)\\ u_t(\cdot,t_1)\end{pmatrix} - \begin{pmatrix} u(\cdot,t_2)\\ u_t(\cdot,t_2)\end{pmatrix}\right\|_{\dot{H}^1 \times L^2(\Rm^5)}\\
\lesssim_1 &\lim_{t_1, t_2 \rightarrow +\infty}  \|F\|_{Z([t_1,t_2])} = 0.
\end{align*}
Therefore there exists $(u_0^+, u_1^+)\in \dot{H}^1 \times L^2$ so that 
\[
 \lim_{t\rightarrow +\infty} \left\|\mathbf{S}_L(-t) \begin{pmatrix} u(\cdot,t)\\ u_t(\cdot,t)\end{pmatrix} - \begin{pmatrix} u_0^+\\ u_1^+\end{pmatrix}\right\|_{\dot{H}^1 \times L^2(\Rm^5)} = 0.
\] 
Thus we have 
\begin{equation} \label{forward scattering}
\lim_{t\rightarrow +\infty} \left\| \begin{pmatrix} u(\cdot,t)\\ u_t(\cdot,t)\end{pmatrix} - \mathbf{S}_L(t) \begin{pmatrix} u_0^+\\ u_1^+\end{pmatrix}\right\|_{\dot{H}^1 \times L^2(\Rm^5)} = 0.
\end{equation}
Let $G^+$ be the radiation profile of the free wave $u_L^+ = \mathbf{S}_L (u_0^+, u_1^+)$ in the positive time direction. By the property of radiation fields we have 
\[
  \lim_{t\rightarrow +\infty} \int_0^\infty \left( \left|G^+(r-t) - r^2 (\partial_t u_L^+) (r, t)\right|^2 + \left|G^+(r-t) + r^2 (\partial_r u_L^+) (r, t)\right|^2 \right) dr  = 0.
\]
We may combine this with \eqref{forward scattering} to conclude 
\[
 \lim_{t\rightarrow +\infty} \int_0^\infty \left(\left|G^+(r-t) - r^2 u_t (r, t)\right|^2 + \left|G^+(r-t) + r^2 u_r (r, t)\right|^2\right) dr = 0.
\]
In addition, if $R>0$ is a constant, then the limit above implies that 
\[
 \|G^+\|_{L^2([R,+\infty))} = \|G^+(r-t)\|_{L_r^2([t+R, +\infty))} = \lim_{t\rightarrow +\infty} \|r^2 u_t(r,t)\|_{L_r^2([t+R, +\infty))}. 
\]
Finally we combine the Strichartz estimates with finite speed of propagation to conclude 
\[
  \|G^+\|_{L^2([R,+\infty))}  \lesssim_1 \lim_{t\rightarrow +\infty} \|u_t(x,t)\|_{L^2(\{x: |x|>t+R\})} 
   \lesssim_1 \lim_{t\rightarrow +\infty} \|\chi_R F\|_{Z([0,t))} = \|\chi_R F\|_{Z([0,+\infty))}. 
\]
The proof in the negative time direction is similar.
\end{proof}

\begin{remark} \label{difference scatter profile}
 The functions $G^+$, $G^-$ are unique. Because they have to be the radiation profiles of the free waves $u_L^\pm$. We call $G^\pm$ the radiation profiles of $u$. If $\tilde{u}$ solves a similar equation 
\[
  \left\{\begin{array}{ll} \partial_t^2 \tilde{u} - \Delta \tilde{u} = \tilde{F}(t,x); & (x,t)\in \Rm^5 \times \Rm; \\
  (\tilde{u},\tilde{u}_t)|_{t=0} = (0,0). & \end{array} \right.
 \]
 with $\|\tilde{F}\|_{Z(\Rm)} < +\infty$, then we may consider its corresponding radiation profiles $\tilde{G}^+$, $\tilde{G}^-$. By linearity $w = u-\tilde{u}$ solves the equation $\partial_t^2 w - \Delta w = F - \tilde{F}$ with zero initial data, whose corresponding radiation profiles are exactly $G^\pm - \tilde{G}^\pm$. We may apply Lemma \ref{scatter profile of nonlinear part} again and obtain 
  \begin{align*}
  \|G^- - \tilde{G}^-\|_{L^2([R,+\infty))} &\lesssim_1 \|\chi_R F - \chi_R \tilde{F}\|_{Z((-\infty,0])};& \\
   \|G^+ - \tilde{G}^+\|_{L^2([R,+\infty))} &\lesssim_1 \|\chi_R F - \chi_R \tilde{F}\|_{Z([0,+\infty))}.&
 \end{align*}
\end{remark}

\begin{lemma}[see, for instance, Li-Shen-Wang-Wei \cite{nonradialCE} for a  proof] \label{recursion lemma}
Assume that $l>1$ and $\alpha>0$ are constants. Let $S: [R,+\infty) \rightarrow [0,+\infty)$ be a function satisfying 
\begin{itemize}
 \item $S(r)\rightarrow 0$ as $r\rightarrow +\infty$;
 \item The recursion formula $S(r_2) \lesssim (r_1/r_2)^\alpha + S^l (r_1)$ holds when $r_2 \gg r_1 \gg R$.
\end{itemize}
Then given any constant $\beta \in (0, (1-1/l)\alpha)$, the decay estimate $S(r) \leq r^{-\beta}$ holds as long as $r>R_0$ is sufficiently large. 
\end{lemma}

\begin{lemma} \label{comparison 1}
 Let $u$ and $v$ be two $R$-weakly non-radiative solutions to (CP1). Then the radiation profiles $G$ and $\tilde{G}$ associated to their initial data satisfy the inequality 
\[
 \|G-\tilde{G}\|_{L^2(\{s: |s|>R\}\times \mathbb{S}^4)} \lesssim_1 \|\chi_R F(t,x,u) - \chi_R F(t,x,v)\|_{Z(\Rm)}. 
\]
\end{lemma}
\begin{proof}
 By the Strichartz estimates we have 
 \begin{align*}
  \sqrt{2} \|G-\tilde{G}\|_{L^2((R,+\infty)\times \mathbb{S}^4)} & = \lim_{t\rightarrow -\infty} \|(u_L, \partial_t u_L) - (v_L, \partial_t v_L)\|_{\dot{H}^1\times L^2(\{x: |x|>R+|t|\})} \\
  & \lesssim_1 \lim_{t\rightarrow -\infty} \|(u, \partial_t u)-(v, \partial_t v)\|_{\dot{H}^1\times L^2(\{x: |x|>R+|t|\})} \\
  & \qquad + \lim_{t\rightarrow -\infty} \|\chi_R F(t,x,u) - \chi_R F(t,x,v)\|_{Z((t,0])}\\
  & \lesssim_1 \|\chi_R F(t,x,u) - \chi_R F(t,x,v)\|_{Z((-\infty,0])}.
 \end{align*}
 Similarly we use the relationship between radiation profiles in two time directions and obtain 
\[
 \|G-\tilde{G}\|_{L^2((R,+\infty)\times \mathbb{S}^4)}^2 \lesssim_1 \|\chi_R F(t,x,u) - \chi_R F(t,x,v)\|_{Z([0,+\infty))}.
\]
Combining these two inequalities, we finish the proof. 
\end{proof}
\begin{lemma} \label{comparison 2} 
Let $u$ and $v$ be two radial $R$-weakly non-radiative solutions to (CP1) with initial data $(u_0,u_1)$ and $(v_0,v_1)$, respectively. If $\|\chi_R u\|_{Y(\Rm)}$, $\|\chi_R v\|_{Y(\Rm)}$ are both sufficiently small, then the radiation profiles $G$ and $\tilde{G}$ of these initial data satisfy 
\[
 \|G-\tilde{G}\|_{L^2(\{s: |s|>R\})} \lesssim_\gamma \left(\|\chi_R u\|_{Y(\Rm)}^{4/3} + \|\chi_R v\|_{Y(\Rm)}^{4/3}\right) \|\chi_R \mathbf{S}_L (u_0-v_0, u_1-v_1)\|_{Y(\Rm)}. 
\]
In addition, if $w_L$ is the free wave whose radiation profile $g$ is given by
\[
 g(s) = \left\{\begin{array}{ll} G(s) - \tilde{G}(s), & |s| \leq R; \\ 0, & |s|>R. \end{array}\right.
\] 
Then we also have 
 \[
 \|G-\tilde{G}\|_{L^2(\{s: |s|>R\}\times \mathbb{S}^4)} \lesssim_\gamma \left(\|\chi_R u\|_{Y(\Rm)}^{4/3} + \|\chi_R v\|_{Y(\Rm)}^{4/3}\right) \|\chi_R w_L\|_{Y(\Rm)}. 
\]
\end{lemma}
\begin{proof}
First of all, the Strichartz estimates gives ($u_L = \mathbf{S}_L (u_0,u_1)$, $v_L = \mathbf{S}_L (v_0,v_1)$)
 \begin{align*}
  \|\chi_R (u-v)\|_{Y(\Rm)}  & \leq \|\chi_R (u_L-v_L)\|_{Y(\Rm)} + C_1 \|\chi_R F(t,x,u) - \chi_R F(t,x,v)\|_{Z(\Rm)}\\
  & \leq \|\chi_R (u_L-v_L)\|_{Y(\Rm)} + C_2 \left(\|\chi_R u\|_{Y(\Rm)}^{4/3} + \|\chi_R u\|_{Y(\Rm)}^{4/3}\right) \|\chi_R (u-v)\|_{Y(\Rm)}
 \end{align*}
 Here the constants $C_1$, $C_2$ depend on the norm $\gamma$ of $F$ only. Thus if $\|\chi_R u\|_{Y(\Rm)}$, $\|\chi_R v\|_{Y(\Rm)}$ are sufficiently small, we always have 
 \[
  \|\chi_R (u-v)\|_{Y(\Rm)} \leq 2 \|\chi_R (u_L-v_L)\|_{Y(\Rm)}. 
 \]
 We then apply Lemma \ref{comparison 1} to obtain the first inequality in the conclusion. 
 \begin{align*}
  \|G-\tilde{G}\|_{L^2(\{s: |s|>R\}\times \mathbb{S}^4)} & \lesssim_1 \|\chi_R F(t,x,u) - \chi_R F(t,x,v)\|_{Z(\Rm)}\\
  & \lesssim_\gamma \left(\|\chi_R u\|_{Y(\Rm)}^{4/3} + \|\chi_R u\|_{Y(\Rm)}^{4/3}\right) \|\chi_R (u-v)\|_{Y(\Rm)}\\
  & \lesssim_\gamma \left(\|\chi_R u\|_{Y(\Rm)}^{4/3} + \|\chi_R u\|_{Y(\Rm)}^{4/3}\right) \|\chi_R (u_L-v_L)\|_{Y(\Rm)}.
 \end{align*}
 In order to prove the second inequality, we observe that if $\tilde{w}$ is the free wave with radiation profile $G - \tilde{G} - g$, then we have $u_L -v_L = w + \tilde{w}$, thus 
 \[
  \|\chi_R (u_L-v_L)\|_{Y(\Rm)} \leq \|\chi_R w\|_{Y(\Rm)} + \|\chi_R \tilde{w}\|_{Y(\Rm)} \leq \|\chi_R w\|_{Y(\Rm)} + \|\tilde{w}\|_{Y(\Rm)}
 \]
 The Strichartz estimates then give
 \[
  \|\tilde{w}\|_{Y(\Rm)} \lesssim_1 \|G - \tilde{G} - g\|_{L^2(\Rm)} \lesssim_1 \|G-\tilde{G}\|_{L^2 (\{s:|s|>R\})}. 
 \]
 Thus we have 
 \[
   \|\chi_R (u_L-v_L)\|_{Y(\Rm)} \lesssim_1 \|\chi_R w\|_{Y(\Rm)} + \|G-\tilde{G}\|_{L^2 (\{s:|s|>R\})}. 
 \]
 Plugging this in the first inequality in the conclusion, we have 
 \[
   \|G-\tilde{G}\|_{L^2(\{s: |s|>R\})} \lesssim_\gamma \left(\|\chi_R u\|_{Y(\Rm)}^{4/3} + \|\chi_R v\|_{Y(\Rm)}^{4/3}\right) \left(\|\chi_R w\|_{Y(\Rm)} + \|G-\tilde{G}\|_{L^2 (\{s:|s|>R\})}\right).
 \]
Our assumption that $\|\chi_R u\|_{Y(\Rm)}$, $\|\chi_R v\|_{Y(\Rm)}$ are both small guarantees that the term 
\[
 \left(\|\chi_R u\|_{Y(\Rm)}^{4/3} + \|\chi_R v\|_{Y(\Rm)}^{4/3}\right) \|G-\tilde{G}\|_{L^2 (\{s:|s|>R\})}
\]
in the right hand side can be absorbed by the left hand side. This finishes the proof. 
 \end{proof}
 
\section{The First Characteristic Number} 
In this section we prove part (a) of our main theorem. We start by 
 \begin{proposition} \label{decay of G}
  Assume that $u$ is a radial $R$-weakly non-radiative solution to (CP1). Let $G$ be the radiation profile associated to the initial data $(u_0,u_1)$. Then for sufficiently large $r\gg R$ we have 
 \begin{align*}
  &\|\chi_r u\|_{Y(\Rm)} \lesssim r^{-1/2}; & &\|G\|_{L^2(\{s: |s|>r\})} \lesssim r^{-7/6}.& 
 \end{align*}
 \end{proposition}
 \begin{proof}
  The proof is similar to the decay estimates of radiation profiles $G$ in 3-dimensional case, as given in \cite{nonradialCE}. Given $r_1\geq R$, we split the radiation profile into two parts
 \begin{align}
   &G_1(s) = \left\{\begin{array}{ll} G(s), & |s|\leq r_1; \\ 0, & |s| > r_1; \end{array}\right. &  &G_2(s) = \left\{\begin{array}{ll} 0, & |s|\leq r_1; \\ G(s), & |s| > r_1; \end{array}\right. & \label{split G}
  \end{align}
 and split the linear free wave $u_L = \mathbf{S}_L (u_0,u_1)$ accordingly: $u_L = u_{1,L} + u_{2,L}$. By Corollary \ref{outward decay} and the Strichartz estimates $S(r) \doteq \|\chi_r u_L\|_{Y(\Rm)}$ satisfies the inequality ($r_2 \geq r_1$)
 \[
  S(r_2) \leq \|\chi_{r_2} u_{1,L}\|_{Y(\Rm)} + \|\chi_{r_2} u_{2,L}\|_{Y(\Rm)} \lesssim_1 (r_1/r_2)^{1/2} \|G_1\|_{L^2(\Rm)} + \|G_2\|_{L^2(\Rm)}
 \]
 Therefore 
 \begin{equation} \label{recurrence 1}
  S(r_2) \lesssim_1 (r_1/r_2)^{1/2} \|G\|_{L^2(\Rm)} + \|G\|_{L^2(\{s: |s|>r_1\})}, \qquad r_2 \geq r_1 \geq R. 
 \end{equation} 
 It immediately follows that $S(r) \rightarrow 0$ as $r \rightarrow +\infty$. If $S(r)$ is sufficiently small, then a combination of Lemma \ref{Y linear to nonlinear} and uniqueness of exterior solutions implies 
 \begin{equation} \label{nonlinear by linear}
  \|\chi_r u\|_{Y(\Rm)} \leq 2S(r). 
 \end{equation}
 Next we assume $r \gg R$ so that $S(r) \ll 1$ is sufficiently small. We apply Lemma  \ref{comparison 1} on $u$ and $v = 0$
 \begin{equation} \label{upper bound G on S}
  \|G\|_{L^2(\{s: |s|>r\})} \lesssim_1 \|\chi_r F(t,x,u)\|_{Z(\Rm)} \lesssim_\gamma \|\chi_r u\|^{7/3} \lesssim_\gamma S^{7/3}(r).
 \end{equation}
 Therefore we plug this upper bound in \eqref{recurrence 1} and obtain the recurrence formula 
 \[
  S(r_2) \lesssim_\gamma (r_1/r_2)^{-1/2} \|G\|_{L^2(\Rm)} + S^{7/3} (r_1), \qquad r_2 \geq r_1 \gg R.  
 \]
 We then apply Lemma \ref{recursion lemma} and obtain a decay estimate 
 \[
  S(r) < r^{-1/4}, \qquad r\gg R. 
 \]
 We then plug this in \eqref{upper bound G on S} to obtain
 \[
  \|G\|_{L^2(\{s: |s|>r\})} \lesssim r^{-7/12}, \qquad r\gg R. 
 \]
 This immediately gives another two estimates ($r\gg R$)
 \begin{align*}
 &\left|\int_{-r}^r sG(s) ds\right| \lesssim r^{11/12};& &\left|\int_{-r}^r G(s) ds\right| \lesssim 1.&
\end{align*}
We then apply Corollary \ref{outward decay 2} and recall the estimates of $G$ given above 
\[
 S(r) \lesssim r^{-3/2} \left|\int_{-r}^r sG(s) ds\right| + r^{-1/2} \left|\int_{-r}^r G(s) ds\right| + \|G\|_{L^2(\{s: |s|>r\})}\lesssim  r^{-1/2}. 
\]
We then finish the proof by \eqref{nonlinear by linear} and \eqref{upper bound G on S}. 
 \end{proof}

\paragraph{The first characteristic number} It immediately follows that if $u$ is a radial $R$-weakly non-radiative solution to (CP1), then the radiation profile $G$ associated to its initial data satisfies $G \in L^1 (\Rm)$. We may define the value of 
 \[
  \alpha = - \int_{-\infty}^{+\infty} G(s) ds 
 \]
 to be the first characteristic number of this solution. Given $r \gg R$, we split $G$ as at the beginning of the proof above at the radius $r$, then write $(u_0,u_1) = (u_0^1, u_1^1) + (u_0^2, u_1^2)$ accordingly. We have
 \[
  \|G_2\|_{L^2(\Rm)} \lesssim r^{-7/6}, \qquad \Rightarrow \qquad \|(u_0^2, u_1^2)\|_{\dot{H}^1\times L^2(\Rm^5)} \lesssim r^{-7/6}.  
 \]
Next we apply Lemma \ref{u by G} and obtain 
\begin{align*}
 &u_0^1(r') = \frac{1}{r'^3} \int_{-r}^r sG(s) ds;& &u_1^1 (r') = - \frac{1}{r'^3} \int_{-r}^r G(s) ds = \frac{\alpha}{r'^3} + \frac{1}{r'^3} \int_{|s|>r} G(s)ds.&
\end{align*}
Thus we have 
\[
 \|(u_0,u_1) - (0, \alpha |x|^{-3})\|_{\dot{H}^1 \times L^2(\{x: |x|>r\})} \lesssim r^{-7/6}.
\]

\begin{lemma} \label{existence of solution with first alpha}
Given any real number $\alpha$, there exists a radial weakly non-radiative solution to (CP1) whose first characteristic number is exactly $\alpha$. 
\end{lemma}
\begin{proof}
Since zero solution is clearly a non-radiative solution with first characteristic number $0$, we assume $\alpha \neq 0$. We consider a complete distance space
\[
 X = \left\{G \in L^2(\{s: |s|>R\}): \|G\|_{L^2(\{s: |s|>r\})} \leq c |\alpha|^{7/3} r^{-7/6}, \forall r \geq R\right\},
\]
whose distance is defined by 
\[
 d(G_1,G_2) = \sup_{r\geq R} r^{7/6} \|G_1-G_2\|_{L^2(\{s: |s|>r\})}. 
\]
Here $c = c(\gamma) \gg 1$ and $R = c^{3/2} |\alpha|^2$ are both positive constants. We also define a map $\mathbf{T}$ from $X$ to itself. Given $G\in X$, we first extend its domain to all $s\in \Rm$ so that 
\[
 G(s) = \frac{-\alpha}{2R} - \frac{1}{2R}\int_{|s'|>R} G(s') ds', \qquad |s|\leq R;.
\]
and let $(u_0,u_1)\in \dot{H}^1 \times L^2(\Rm^5)$ be initial data so that the corresponding radiation profile is exactly $G$. Let $u$ be the solution to the non-linear wave equation 
\[ 
 \partial_t^2 u - \Delta u = \chi_R F(t,x,u)
\]
By Corollary \ref{outward decay 2}, the free wave $u_L = \mathbf{S}_L (u_0,u_1)$ satisfies ($r\geq R$)
\begin{align*}
 \|\chi_r u_L\|_{Y(\Rm)} & \lesssim_1 r^{-3/2} \left|\int_{-r}^r sG(s) ds\right| +  r^{-1/2} \left|\int_{-r}^r G(s) ds\right| +\|G\|_{L^2(\{s: |s|>r\})}\\
 & \lesssim_1 r^{-3/2} \cdot c|\alpha|^{7/3} r^{1/3} + |\alpha| r^{-1/2} + r^{-1/2}  \left|\int_{|s|>r} G(s) ds\right| + c |\alpha|^{7/3} r^{-7/6}\\
 & \lesssim_1 |\alpha| r^{-1/2} (1 + c|\alpha|^{4/3} r^{-2/3})  \\
 & \lesssim_1 |\alpha| r^{-1/2}.
\end{align*}
In particular, $\|\chi_R u_L\|_{Y(\Rm)} \lesssim_1 |\alpha| R^{-1/2} \ll 1$. A combination of Lemma \ref{Y linear to nonlinear} and the global extension method then guarantee that $u$ is a globally defined solution satisfying $\|\chi_r u\|_{Y(\Rm)} \lesssim_1 |\alpha| r^{-1/2}$ for all $r\geq R$. Next we consider $u_{NL} = u - u_L$ and apply Lemma \ref{scatter profile of nonlinear part}. This immediately gives two radiation profiles $G^\pm$ so that 
 \begin{align*}
  \lim_{t\rightarrow +\infty} \int_0^\infty \left(\left|G^+(r-t) - r^2 \partial_t u_{NL} (r, t)\right|^2 + \left|G^+(r-t) + r^2 \partial_r u_{NL} (r, t)\right|^2\right) dr & = 0; \\
  \lim_{t\rightarrow -\infty} \int_0^\infty \left(\left|G^-(t+r) - r^2 \partial_t u_{NL} (r,t)\right|^2 +  \left|G^-(t+r) - r^2 \partial_r u_{NL}(r,t)\right|^2\right)dr & = 0.
 \end{align*}
In addition we have 
\begin{align*}
  \|G^-\|_{L^2((r,+\infty))} &\lesssim_1 \|\chi_r F(t,x,u)\|_{Z((-\infty, 0])} \lesssim_\gamma \|\chi_r u\|_{Y(\Rm)}^{7/3} \lesssim_\gamma |\alpha|^{7/3} r^{-7/6}, & & r\geq R;\\
  \|G^+\|_{L^2((r,+\infty))} &\lesssim_1 \|\chi_r F(t,x,u)\|_{Z([0,+\infty))} \lesssim_\gamma \|\chi_r u\|_{Y(\Rm)}^{7/3} \lesssim_\gamma |\alpha|^{7/3} r^{-7/6}, & & r\geq R.
\end{align*}
We define $\mathbf{T} : X \rightarrow X$ by 
\[
 (\mathbf{T} G) (s) = \left\{\begin{array}{ll} -G^-(s), & s>R; \\ -G^+(-s), & s<-R. \end{array}\right.
\]
This is still in the space $X$ by our assumption $c \gg 1$ and the upper bounds of $G^\pm$ given above. Next we verify this map $\mathbf{T}$ is a contraction map. Given another element $\tilde{G} \in X$, we define $\tilde{u}, \tilde{u}_L, \tilde{u}_{NL}, \tilde{G}^\pm$ as above. We may apply Corollary \ref{outward decay 2} and obtain ($r\geq R$)
\begin{align*}
 \left\|\chi_r (u_L- \tilde{u}_L) \right\|_{Y(\Rm)} & \lesssim_1  r^{-3/2} \left|\int_{-r}^r s(G-\tilde{G}) ds\right| +  r^{-1/2} \left|\int_{-r}^r (G-\tilde{G}) ds\right| +\|G-\tilde{G}\|_{L^2(\{s: |s|>r\})}\\
 & \lesssim_1  r^{-3/2} \left|\int_{-r}^r s(G-\tilde{G}) ds\right| +  r^{-1/2} \left|\int_{|s|>r} (G-\tilde{G}) ds\right| +\|G-\tilde{G}\|_{L^2(\{s: |s|>r\})} \\
 & \lesssim_1 r^{-7/6} d(G, \tilde{G}).
\end{align*}
We then apply Lemma \ref{Y linear to nonlinear} and the uniqueness of exterior solutions to obtain 
\[
 \|\chi_r u - \chi_r \tilde{u}\|_{Y(\Rm)} \lesssim_1 r^{-7/6} d(G,\tilde{G}), \qquad r\geq R.
\]
Remark \ref{difference scatter profile} implies that if $r\geq R$, then 
\begin{align*}
 \|G^\pm - \tilde{G}^\pm\|_{L^2((r,+\infty))} & \lesssim_1 \|\chi_r F(t,x,u) - \chi_r F(t,x,\tilde{u})\|_{Z(\Rm)}\\
 & \lesssim_\gamma \left(\|\chi_r u\|_{Y(\Rm)}^{4/3} +  \|\chi_r \tilde{u}\|_{Y(\Rm)}^{4/3}\right)\|\chi_r u - \chi_r \tilde{u}\|_{Y(\Rm)}\\
 & \lesssim_\gamma |\alpha|^{4/3} r^{-11/6} d(G,\tilde{G}). 
\end{align*}
Therefore 
\[
 d(\mathbf{T} G, \mathbf{T} \tilde{G}) \lesssim_\gamma   |\alpha|^{4/3} R^{-2/3} d(G,\tilde{G}).
\]
Our assumption on $c$ and $R$ then guarantees that $\mathbf{T}$ is a contraction map from $X$ to itself. Therefore there exists a fixed-point $G$ so that 
$\mathbf{T} G = G$. We claim that the corresponding solution $u$ to (CP1) is an $R$-weakly non-radiative solution. In fact we have 
\begin{align*}
 \lim_{t\rightarrow \pm\infty} \|\nabla_{x,t} u(\cdot,t)\|_{L^2(\{x: |x|>R+|t|\})}  = \lim_{t\rightarrow \pm\infty} \|\nabla_{x,t} (u_L (\cdot,t)+u_{NL}(\cdot,t))\|_{L^2(\{x: |x|>R+|t|\})}.
\end{align*} 
We recall the basic theory of radiation fields
\begin{align*}
 \lim_{t\rightarrow +\infty} \|\nabla_{x,t} u(\cdot,t)\|_{L^2(\{x: |x|>R+|t|\})} & = \sqrt{2} \sigma_4^{1/2} \|G(-s)+G^+(s)\|_{L^2((+R,+\infty))} = 0; \\
 \lim_{t\rightarrow -\infty} \|\nabla_{x,t} u(\cdot,t)\|_{L^2(\{x: |x|>R+|t|\})} & = \sqrt{2} \sigma_4^{1/2} \|G+G^-\|_{L^2((+R, +\infty))} = 0. 
 \end{align*}
 Finally the way in which we define $G$ in $[-R,R]$ guarantees that the first characteristic number of $u$ is exactly $\alpha$. 
\end{proof}

\begin{remark} 
 If $F(t,x,u)$ is even in $t$ and odd in $u$, i.e. satisfies the assumption (AS), then the weakly non-radiative solutions $u$ constructed in Lemma \ref{existence of solution with first alpha} always satisfy
 \begin{itemize}
  \item The solutions are odd in time $t$, i.e. $u(x,-t) = -u(x,t)$;
  \item The radiation profiles of their initial data are even, i.e. $G(-s) = G(s)$.  
 \end{itemize}
\end{remark}
\begin{proof}
We first verify that if $G$ is even, then the corresponding solution $u$ defined in the proof above is odd in time. In fact, if $G$ is even, then by Lemma \ref{u by G} the initial data satisfy $u_0 = 0$. As a result, the function $-u(x,-t)$ is also a solution to the same equation with the same initial data
\[
 (\partial_t^2 - \Delta) (-u(x,-t)) = - (\partial_t^2 u - \Delta u) (x,-t) = -\chi_R F(-t,x,u(x,-t)) = \chi_R F(t,x,-u(x,-t)). 
\]
Thus we must have $u(x,-t) = -u(x,t)$ by the uniqueness of the solution. As a result, it suffices to prove that $G$ is even. 
 Because the radiation profile $G$ is constructed via a contraction map, $G$ can be obtained by a limit process in the space $X$
 \[
  G = \lim_{k\rightarrow +\infty} G_k.
 \]
 Here the sequence $G_k$ can be defined by 
 \begin{align*}
  &G_0 = 0;& &G_{k+1} = \mathbf{T} G_k.&
 \end{align*} 
 We only need to prove each $G_k$ is an even function of $s$ by an induction. Let us assume $G$ is even and show that $\mathbf{T}G$ is also even. We have already shown that the corresponding solution $u$ must be odd in time. This also means $u_{NL}$ is odd in time. By the definition of radiation fields we have the limit in the space $L^2([R,+\infty))$:
 \[
  G^- =  \lim_{t\rightarrow -\infty} (s-t)^2 (\partial_t u_{NL}) (s-t, t)
 \]
Since $u_{NL}$ is odd, $\partial_t u_{NL}$ is even. Thus  
\begin{align*}
 G^-  = \lim_{t\rightarrow +\infty} (s+t)^2 (\partial_t u_{NL}) (s+t,-t) = \lim_{t\rightarrow +\infty} (s+t)^2 (\partial_t u_{NL}) (s+t,t) = G^+ (s).
\end{align*}
Therefore $\mathbf{T} G$ is still even. 
\end{proof}
\section{The Second Characteristic Number}

In this section we prove parts (b) and (c) of our main theorem. 

\begin{proposition} \label{comparison with same first characteristic}
 Let $u$ and $v$ be two radial $R$-weakly non-radiative solutions to (CP1) with the same first characteristic number. Then the radiation profiles $G$ and $\tilde{G}$ of their initial data satisfy the inequality 
\[
 \|G-\tilde{G}\|_{L^2(\{s: |s|>R\})} \lesssim r^{-13/6}, \qquad r\gg R. 
\]
\end{proposition}
\begin{proof}
 We apply Lemma \ref{comparison 2} for sufficiently large $r \gg R$. 
 \begin{equation} \label{difference of G 1}
 \|G-\tilde{G}\|_{L^2(\{s: |s|>r\})} \lesssim \left(\|\chi_r u\|_{Y(\Rm)}^{4/3} + \|\chi_r v\|_{Y(\Rm)}^{4/3}\right) \|\chi_r w_L\|_{Y(\Rm)}. 
\end{equation}
 Here $w_L$ is the linear free wave with radiation profile 
 \[
 g(s) = \left\{\begin{array}{ll} G(s) - \tilde{G}(s), & |s| \leq r; \\ 0, & |s|>r. \end{array}\right.
\] 
Lemma \ref{outward decay} immediately gives a decay estimate 
\begin{align}
 \|\chi_r w_L\|_{Y(\Rm)} & \lesssim r^{-3/2} \left|\int_{-r}^r s[G(s)-\tilde{G}(s)] ds\right| + r^{-1/2} \left|\int_{-r}^r [G(s)-\tilde{G}(s)] ds\right| \nonumber\\
 & \lesssim r^{-3/2} \left|\int_{-r}^r s[G(s)-\tilde{G}(s)] ds\right| + r^{-1/2}  \left|\int_{|s|>r} [G(s)-\tilde{G}(s)] ds\right|. \label{upper bound of wL 1}
\end{align} 
Here we use the assumption on the first characteristic numbers. Next we recall the decay of $G$, $\tilde{G}$ and obtain $\|\chi_r w_L\|_{Y(\Rm)} \lesssim  r^{-7/6}$ for sufficiently large $r$. This immediately gives a decay estimate of $G - \tilde{G}$ by \eqref{difference of G 1} and the decay estimates $\|\chi_r u\|_{Y(\Rm)}, \|\chi_r v\|_{Y(\Rm)} \lesssim r^{-1/2}$ given in Proposition \ref{decay of G}.
\[
 \|G-\tilde{G}\|_{L^2(\{s: |s|>r\})} \lesssim r^{-11/6}, \qquad r\gg R. 
\]
This enable us to give a better estimate of $w_L$ by \eqref{upper bound of wL 1}
\[
 \|\chi_r w_L\|_{Y(\Rm)} \lesssim r^{-3/2}, \qquad r\gg R. 
\]
Finally we use \eqref{difference of G 1} again to finish the proof. 
\end{proof}

\paragraph{The second characteristic number} Let us fix $\alpha$ and a radial weakly non-radiative solution $u^\alpha$. Let $G_\alpha$ be the radiation profile to the initial data of $u^\alpha$. If $u$ is another radial weakly non-radiative solution to (CP1) with the same first characteristic number $\alpha$, then the radiation profile $G$ to its initial data $(u_0,u_1)$ satisfies 
 \begin{align*} 
   &\|G-G_\alpha\|_{L^2(\{s: |s|>r\})} \lesssim r^{-13/6}, \quad r\gg 1;& &\Rightarrow& &sG(s)-s\tilde{G}(s) \in L^1 (\Rm).&
  \end{align*}
We call the number 
\[
 \beta = \int_{-\infty}^{+\infty} [sG(s)-s\tilde{G}(s)] ds
\]
the second characteristic number of $u$ (with respect to $u^\alpha$). 

\paragraph{Asymptotic behaviour of initial data} We then describe the asymptotic behaviour of $(u_0,u_1)$ in term of the second characteristic number $\beta$. Given $r\gg 1$, we first choose $G_{0,\beta}$ be radiation profile so that $G_{0,\beta}(s) = 0$ for $|s|>r$ and 
\begin{align*}
 &\int_{-r}^r G_{0,\beta}(s) ds = 0;& &\int_{-r}^r s G_{0,\beta} (s) ds = \beta.& 
\end{align*}
By Lemma \ref{u by G} we have the corresponding initial data $(u_{0,\beta} , u_{1,\beta})$ satisfies 
\begin{equation} \label{u0u1 outward region}
 (u_{0,\beta}(x), u_{1,\beta}(x)) = (\beta |x|^{-3}, 0), \qquad |x|>r.
\end{equation}
Thus we have (for simplicity we use the notation $\mathcal{H}(r) = \dot{H}^1 \times L^2(\{x: |x|>r\})$)
\[
 \|(u_0,u_1) - (u_0^\alpha, u_1^\alpha) - (\beta |x|^{-3},0)\|_{\mathcal{H}(r)} = \|(u_0,u_1) - (u_0^\alpha, u_1^\alpha) - (u_{0,\beta}(x), u_{1,\beta}(x))\|_{\mathcal{H}(r)}
\]
The radiation profile of $(u_0,u_1) - (u_0^\alpha, u_1^\alpha) - (u_{0,\beta} , u_{1,\beta})$ is $G - G_\alpha - G_{0,\beta}$. We then apply Corollary \ref{outward decay 2} to give an upper bound of the $\mathcal{H}_r$ norm above (up to an absolute constant)
\[
 r^{-3/2} \left|\int_{-r}^r s (G-G_\alpha - G_{0,\beta}) ds \right| + r^{-1/2} \left|\int_{-r}^r (G-G_\alpha-G_{0,\beta}) ds\right| + \|G-G_\alpha - G_{0,\beta}\|_{L^2(\{s: |s|>r\})}.
\]
The assumptions on the characteristic numbers imply
\begin{align*}
 &\int_{-\infty}^{+\infty} s(G-G_\alpha) ds = \beta;& &\int_{-\infty}^{+\infty} (G-G_\alpha) ds = 0.& 
\end{align*}
We combine these identities with the assumptions on $G_{0,\beta}$ to rewrite the upper bound above as below
\[
 r^{-3/2} \left|\int_{|s|>r} s (G-G_\alpha) ds \right| + r^{-1/2} \left|\int_{|s|>r} (G-G_\alpha) ds\right| + \|G-G_\alpha\|_{L^2(\{s: |s|>r\})}.
\]
Finally we utilize the decay estimate $\|G-G_\alpha\|_{L^2(\{s: |s|>r\})} \lesssim r^{-13/6}$ to conclude that  
\[
 \|(u_0,u_1) - (u_0^\alpha, u_1^\alpha) - (\beta |x|^{-3},0)\|_{\mathcal{H}(r)} \lesssim r^{-13/6}. 
\]
as $r$ is sufficiently large. 

\begin{proposition} \label{existence of solution with second beta}
 Given two real numbers $\alpha$, $\beta$, let $u^\alpha$ be the weakly non-radiative solution with first characteristic number $\alpha$ given in Lemma \ref{existence of solution with first alpha}. Then there exists a radial weakly non-radiative solution $u$ whose second characteristic number is exactly $\beta$. 
\end{proposition} 
\begin{proof}
 The proof is similar to Lemma \ref{existence of solution with first alpha}. We also use the notation $G_\alpha$ for radiation profile of the initial data $(u_0^\alpha, u_1^\alpha)$ of $u^\alpha$. We recall that the linear free wave $u_L^\alpha = \mathbf{S_L}(u_0^\alpha, u_1^\alpha)$ and $u^\alpha$ satisfies 
 \begin{align*}
  &\|\chi_r u_L^\alpha\|_{Y(\Rm)} \lesssim_1  |\alpha| r^{-1/2},& &\|\chi_r u^\alpha\|_{Y(\Rm)} \lesssim_1 |\alpha| r^{-1/2},& &r \geq R_0.&
 \end{align*}
 Here $R_0 \simeq_\gamma |\alpha|^2$. Since $u^\alpha$ itself is a weakly non-radiative solution with second characteristic number zero, we may assume $\beta \neq 0$. We define a distance space 
 \begin{align*}
  X  = \left\{G\in L^2(\{s: |s|>R\}): \begin{array}{ll}\|G\|_{L^2(\{s: |s|>r\})}\leq c |\beta|^{7/3} r^{-7/2}, & r\geq R, r\leq |\beta|/|\alpha|\\
   \|G\|_{L^2(\{s: |s|>r\})}\leq c|\alpha|^{4/3} |\beta| r^{-13/6}, & r\geq R, r > |\beta|/|\alpha|\end{array}\right\},
 \end{align*}
 whose distance $d(G_1,G_2)$ is defined by 
 \[
   \max\left\{\sup_{r\geq R, r\leq \frac{|\beta|}{|\alpha|}} |\beta|^{-\frac{7}{3}} r^{\frac{7}{2}} \|G_1-G_2\|_{L^2(\{s: |s|>r\})}, \sup_{r\geq R, r\geq \frac{|\beta|}{|\alpha|}} |\alpha|^{-\frac{4}{3}} |\beta|^{-1} r^{\frac{13}{6}}\|G_1-G_2\|_{L^2(\{s: |s|>r\})} \right\}.
 \]
 Here the constants $c = c(\gamma) \gg 1$ is a large constant and $R = \max \{c^{3/2} |\alpha|^2, c^{1/2}|\beta|^{2/3}\}\geq R_0$. If $\alpha = 0$, then we understand $|\beta|/|\alpha| = +\infty$ and ignore the case $r > |\beta|/|\alpha|$ above. Similarly if $|\beta|/|\alpha| \leq R$, we ignore the case $r \leq |\beta|/|\alpha|$. Next we define the map $\mathbf{T}: X \rightarrow X$. We first extend the domain of $G$ to the whole real number in a manner so that 
 \begin{align*}
  &\int_{-\infty}^{+\infty} G(s) ds = 0;& &\int_{-\infty}^{+\infty} sG(s) ds = \beta.&
 \end{align*}
 Let $(u_0,u_1)$ be the initial data with radiation profile $G + G_\alpha$. Then the free wave $w_L = \mathbf{S}_L (u_0-u_0^\alpha, u_1-u_1^\alpha)$, whose radiation profile is exactly $G$, satisfies ($r \geq R$)
\begin{align*}
 \|\chi_r w_L\|_{Y(\Rm)} & \lesssim_1 r^{-3/2} \left|\int_{-r}^r sG(s) ds \right| + r^{-1/2} \left|\int_{-r}^r G(s) ds\right| + \|G\|_{L^2(\{s: |s|>r\})} \\
 & \lesssim_1 |\beta| r^{-3/2} + r^{-3/2} \left|\int_{|s|>r} sG(s) ds \right| + r^{-1/2} \left|\int_{|s|>r} G(s) ds\right| \\
 & \lesssim_1 |\beta| r^{-3/2}. 
\end{align*}
Here we use Lemma \ref{outward decay 2}, the decay of $\|u\|_{L^2(\{s: |s|>r\})}$ given in the definition of $X$ and the lower bound of $R$. Thus if $r\geq R$, then we have
\begin{align*}
 \|\chi_r \mathbf{S}_L (u_0,u_1)\|_{Y(\Rm)} &\leq \|\chi_r w_L\|_{Y(\Rm)} + \|\chi_r u_L^\alpha\|_{Y(\Rm)}  \lesssim_1 |\beta| r^{-3/2} \ll 1, \quad R\leq r \leq |\beta|/|\alpha|;\\
 \|\chi_r \mathbf{S}_L (u_0,u_1)\|_{Y(\Rm)} & \leq \|\chi_r w_L\|_{Y(\Rm)} + \|\chi_r u_L^\alpha\|_{Y(\Rm)} \lesssim_1 |\alpha| r^{-1/2} \ll 1, \quad r \geq \max\{R, |\beta|/|\alpha|\}. 
\end{align*}
We consider the solution $u$ to the non-linear wave equation 
\[ 
 \partial_t^2 u - \Delta u = \chi_R F(t,x,u)
\]
with initial data $(u_0,u_1)$. A combination of Lemma \ref{Y linear to nonlinear}, the global extension method and the uniqueness of exterior solutions gives
\begin{equation} \label{Y norms of u}
  \|\chi_r u\|_{Y(\Rm)} \lesssim_1 \left\{\begin{array}{ll} |\beta| r^{-3/2}, & R\leq r \leq |\beta|/|\alpha|; \\ |\alpha| r^{-1/2}, & r\geq \max\{R, |\beta|/|\alpha|\}. \end{array}\right. 
\end{equation}
Similarly we have
\[
 \|\chi_r (u-u^\alpha)\|_{Y(\Rm)} \leq 2\|\chi_r w_L\|_{Y(\Rm)} \lesssim_1 |\beta| r^{-3/2}, \qquad r\geq R. 
\]
Since $w_{NL} = u - u^\alpha - w_L$ solves the wave equation 
\[
 \partial_t^2 w_{NL} - \Delta w_{NL} = \chi_R F(t,x,u) - \chi_{R_0} F(t,x,u^\alpha)
\]
with zero initial data. According to Lemma \ref{scatter profile of nonlinear part}, there exist $G^+, G^- \in L^2(\Rm)$ so that 
\begin{align}
  \lim_{t\rightarrow +\infty} \int_0^\infty \left(\left|G^+(r-t) - r^2 \partial_t w_{NL} (r, t)\right|^2 + \left|G^+(r-t) + r^2 \partial_r w_{NL} (r, t)\right|^2\right) dr & = 0; \label{G plus property} \\
  \lim_{t\rightarrow -\infty} \int_0^\infty \left(\left|G^-(r+t) - r^2 \partial_t w_{NL}(r,t)\right|^2 +  \left|G^-(r+t) - r^2 \partial_r w_{NL}(r,t)\right|^2\right)dr & = 0. \label{G minus property}
 \end{align}
 The profiles $G^\pm$ also satisfy ($r \geq R$)
 \begin{align*}
  \|G^\pm\|_{L^2([r,+\infty))} & \lesssim_1 \|\chi_r F(t,x,u) - \chi_r F(t,x,u^\alpha)\|_{Z(\Rm)} \\
   & \lesssim_\gamma \left(\|\chi_r u\|_{Y(\Rm)}^{4/3} + \|\chi_r u^\alpha\|_{Y(\Rm)}^{4/3}\right) \|\chi_r u - \chi_r u^\alpha\|_{Y(\Rm)}.
 \end{align*}
 Therefore we have 
 \begin{align*}
  \|G^\pm\|_{L^2([r,+\infty))} & \lesssim_\gamma |\beta|^{7/3} r^{-7/2}, & &R\leq r \leq |\beta|/|\alpha|;\\
   \|G^\pm\|_{L^2([r,+\infty))} & \lesssim_\gamma |\alpha|^{4/3}|\beta| r^{-13/6}, & &r \geq \max\{|\beta|/|\alpha|,R\}.
 \end{align*}
Finally we define 
 \[
  (\mathbf{T} G) = \left\{\begin{array}{ll} -G^-(s), & s>R; \\ -G^+(-s), & s<-R. \end{array} \right. 
 \]
 The upper bounds of $L^2$ norms given above and our assumptions on $c$ guarantees that $\mathbf{T} G$ is still contained in the space $X$. Next we verify $\mathbf{T}$ is a contraction map. Let $\tilde{G} \in X$, we define $\tilde{w}_L$, $\tilde{u}$ and $\tilde{G}^\pm$ accordingly. We may apply Corollary \ref{outward decay 2} again and obtain ($r\geq R$)
 \begin{align*}
  \|\chi_r u - \chi_r \tilde{u}\|_{Y(\Rm)} &\lesssim_1 r^{-3/2} \left|\int_{-r}^r s(G - \tilde{G}) ds \right| + r^{-1/2} \left|\int_{-r}^r (G-\tilde{G}) ds\right| + \|G-\tilde{G}\|_{L^2(\{s: |s|>r\})}\\
 & \lesssim_1 r^{-3/2} \left|\int_{|s|>r} s(G - \tilde{G}) ds \right| + r^{-1/2} \left|\int_{|s|>r} (G-\tilde{G}) ds\right| + \|G-\tilde{G}\|_{L^2(\{s: |s|>r\})}\\
 & \lesssim_1 c^{-1} d(G, \tilde{G}) |\beta| r^{-3/2}.
 \end{align*}
A combination of Lemma \ref{Y linear to nonlinear} and the uniqueness of exterior solutions gives 
\[
 \|\chi_r u - \chi_r \tilde{u}\|_{Y(\Rm)} \lesssim_1 d(G,\tilde{G}) c^{-1} |\beta| r^{-3/2}, \qquad r\geq R.
\]
The function $\tilde{w}_{NL} \doteq \tilde{u} - u^\alpha - \tilde{w}_L$ solves the equation 
\[
 \partial_t^2 \tilde{w}_{NL} - \Delta \tilde{w}_{NL} = \chi_R F(t,x,\tilde{u}) - \chi_{R_0} F(t,x,u^\alpha)
\]
with zero initial data. We then recall Remark \ref{difference scatter profile} and obtain ($r\geq R$)
\begin{align*}
 \|G^+ - \tilde{G}^+\|_{L^2([r,+\infty))} & \lesssim_1 \|\chi_r F(t,x,u) - \chi_r F(t,x,\tilde{u})\|_{Z([0,+\infty))} \\
 & \lesssim_\gamma  \left(\|\chi_r u\|_{Y(\Rm)}^{4/3} + \|\chi_r \tilde{u}\|_{Y(\Rm)}^{4/3}\right) \|\chi_r u - \chi_r \tilde{u}\|_{Y(\Rm)}.
\end{align*}
Thus we recall \eqref{Y norms of u} and obtain
 \begin{align*}
  \|G^+ - \tilde{G}^+\|_{L^2([r,+\infty)} & \lesssim_\gamma d(G,\tilde{G}) c^{-1} |\beta|^{7/3} r^{-7/2}, & &R\leq r \leq |\beta|/|\alpha|;\\
   \|G^+ - \tilde{G}^+\|_{L^2([r,+\infty))} & \lesssim_\gamma d(G,\tilde{G}) c^{-1} |\alpha|^{4/3}|\beta| r^{-13/6}, & &r \geq \max\{|\beta|/|\alpha|,R\}.
 \end{align*}
The case of $G^- - \tilde{G}^-$ is similar. Thus 
\[
 d(\mathbf{T} G,\mathbf{T} \tilde{G}) \lesssim_\gamma c^{-1} d(G, \tilde{G}).
\]
Therefore our choice of $c = c(\gamma) \gg 1$ guarantees that the map $\mathbf{T}$ is a contraction map thus comes with a unique fixed point $G$. Next we use the idenity $u = u^\alpha + w_{NL} + w_L$, recall that $u^\alpha$ is a $R$-weakly non-radiative solution, and obtain
\begin{align*}
 \lim_{t\rightarrow +\infty} \int_{|x|>R+|t|} |\nabla_{t,x} u(x,t)|^2 dx & =  \lim_{t\rightarrow +\infty} \int_{|x|>R+|t|} |\nabla_{t,x} (w_{NL}+w_L)(x,t)|^2 dx.
\end{align*}
We use the polar coordinates
\begin{align*}
 \lim_{t\rightarrow +\infty} \int_{|x|>R+|t|} |\nabla_{t,x} u(x,t)|^2 dx  & = \sigma_4 \lim_{t\rightarrow +\infty} \|r^2 \nabla_{t,r} (w_{NL}+w_L)(r,t)\|_{L^2([t+R,+\infty))}^2.
\end{align*}
We then recall \eqref{G plus property} as well as the fact that the radiation profile of $w_L$ is $G$, and conclude
\begin{align*}
 \lim_{t\rightarrow +\infty} \int_{|x|>R+|t|} |\nabla_{t,x} u(x,t)|^2 dx & = 2\sigma_4 \lim_{t\rightarrow +\infty} \|G^+(r-t)+G(t-r)\|_{L^2([t+R,+\infty))}^2 = 0.
\end{align*} 
Thus $u$ is an $R$-weakly non-radiative solution to (CP1). Its first and second characteristic numbers are clearly $\alpha$ and $\beta$ by our construction method. 
 \end{proof}

\begin{corollary} 
 Let $\alpha, \beta \in \Rm$ and $v^\alpha$ be an arbitrary radial weakly non-radiative solution with first characteristic number $\alpha$. Then there exists a radial weakly non-radiative solution whose first characteristic number is also $\alpha$, and whose second characteristic number with respect to $v^\alpha$ is $\beta$. 
\end{corollary}
\begin{proof}
 Let $u^\alpha$ be the radial weakly non-radiative solution as given in Lemma \ref{existence of solution with first alpha}, If the second characteristic number of $v^\alpha$ with respect to $u^\alpha$ is $\beta_0$, then by Proposition \ref{existence of solution with second beta} we may find a radial weakly non-radiative solution whose second characteristic number with respect to $u^\alpha$ is exactly $\beta+\beta_0$. This is clearly a solution with second characteristic number $\beta$ with respect to $v^\alpha$. 
\end{proof}
\begin{proposition} \label{uniqueness of solution}
Let $u$ and $v$ be two radial $R$-weakly non-radiative solutions to (CP1) with the same first characteristic number. If the radiation profiles $G$ and $\tilde{G}$ associated to their initial data satisfy (i.e. these two solution share the same second characteristic number)
\[
 \int_{-\infty}^{+\infty} s[G(s)-\tilde{G}(s)] ds = 0, 
\]
then 
\begin{align*}
 &G(s) = \tilde{G} (s), \quad \hbox{a.e.} \; |s|>R;& &u(x,t) = v(x,t), \; \hbox{if} \; |x|>R+|t|.&
\end{align*}
\end{proposition}
\begin{proof}
 The proof consists of two steps. 
 \begin{itemize} 
  \item We first show that the two identities above holds if we substitute $R$ by a sufficiently large number $R_1 > R$. 
  \item We then show that if the two identities above hold for $r>R$, then they also hold for some $r' < r$. 
 \end{itemize}
 We start with the first step. We first recall $\|\chi_r u\|_{Y(\Rm)}, \|\chi_r v\|_{Y(\Rm)} \lesssim r^{-1/2}$ for sufficiently large $r$ and apply Lemma \ref{comparison 2}. 
 \begin{equation} \label{difference of G 2}
 \|G-\tilde{G}\|_{L^2(\{s: |s|>r\})} \lesssim r^{-2/3} \|\chi_r w_L\|_{Y(\Rm)}, \qquad r\gg R. 
\end{equation}
Here $w_L$ is the linear free wave with radiation profile $G - \tilde{G}$. Lemma \ref{outward decay 2} immediately gives a decay estimate 
\begin{align*}
 \|\chi_r w_L\|_{Y(\Rm)} & \lesssim_1 r^{-3/2} \left|\int_{-r}^r s(G-\tilde{G}) ds\right| + r^{-1/2} \left|\int_{-r}^r (G-\tilde{G}) ds\right| + \|G-\tilde{G}\|_{L^2(\{s: |s|>r\})}\\
 & \lesssim_1 r^{-3/2} \left|\int_{|s|>r} s(G-\tilde{G}) ds\right| + r^{-1/2}  \left|\int_{|s|>r} (G-\tilde{G})ds\right| + \|G-\tilde{G}\|_{L^2(\{s: |s|>r\})}.
\end{align*} 
This implies that ($\kappa \geq 2$)
\begin{equation} \label{induction 1}
 \|G-\tilde{G}\|_{L^2(\{s: |s|>r\})} \leq r^{-\kappa}, \quad \forall\, r\geq r' \qquad \Rightarrow \qquad \|\chi_r w_L\|_{Y(\Rm)} \leq c_2 r^{-\kappa}, \quad \forall\, r\geq r'. 
\end{equation}
Here the constant $c_2$ is an absolute constant, $r'$ is any positive constant. In addition, the conclusion of Proposition \ref{comparison with same first characteristic} and \eqref{difference of G 2} implies that we may find a large $R_1 \gg \max\{R,1\}$, so that 
\[
 \|G- \tilde{G}\|_{L^2(\{s: |s|>r\})} \leq r^{-2}, \qquad \forall r\geq R_1;
\]
and 
\begin{equation} \label{induction 2}
 \|G-\tilde{G}\|_{L^2(\{s: |s|>r\})} \leq c_2^{-1} r^{-1/2} \|\chi_r w_L\|_{Y(\Rm)}, \qquad \forall r\geq R_1.
\end{equation} 
We may combine \eqref{induction 1} and \eqref{induction 2} to obtain ($\kappa \geq 2$)
\[
 \|G-\tilde{G}\|_{L^2(\{s: |s|>r\})} \leq r^{-\kappa}, \quad \forall\, r\geq R_1 \qquad \Rightarrow \qquad \|G-\tilde{G}\|_{L^2(\{s: |s|>r\})} \leq r^{-\kappa-1/2}, \quad \forall\, r\geq R_1. 
\]
An induction in $\kappa$ immediately gives
\[
 \|G-\tilde{G}\|_{L^2(\{s: |s|>r\})} \leq r^{-\kappa}, \qquad \forall r\geq R_1, \; \kappa \geq 2. 
\]
We can make $\gamma \rightarrow +\infty$ and conclude $\|G-\tilde{G}\|_{L^2(\{s: |s|>R_1\})} = 0$. We then apply Lemma \ref{u by G} and obtain 
$(u_0(x),u_1(x)) = (v_0(x), v_1(x))$ if $|x|>R_1$. The finite speed of propagation then gives $u(x,t) = v(x,t)$ if $|x|>R_1+|t|$. This finishes the first step. Next we assume ($r>R$)
\begin{align*}
 &G(s) = \tilde{G} (s), \quad \hbox{a.e.} \; |s|>r;& &u(x,t) = v(x,t), \; \hbox{if} \; |x|>r+|t|;&
\end{align*}
and prove that there exists a radius $r' \in [R,r)$ so that these identities also hold if we substitute $r$ by $r'$. We first apply Lemma \ref{comparison 1}. 
\begin{equation} \label{coincidence G 1}
 \|G-\tilde{G}\|_{L^2(\{s: |s|>r'\})} \lesssim_1 \|\chi_{r'} F(t,x,u) - \chi_{r'} F(t,x,v)\|_{Z(\Rm)}.
\end{equation}
Since $u$ and $v$ coincides in the exterior region $\{x: |x|>r+|t|\}$, we have 
\begin{align}
 \|\chi_{r'} F(t,x,u) - \chi_{r'} F(t,x,v)\|_{Z(\Rm)} & = \|\chi_{r', r} F(t,x,u) - \chi_{r',r} F(t,x,v)\|_{Z(\Rm)} \nonumber\\
 & \lesssim_\gamma \left(\|\chi_{r',r} u\|_{Y(\Rm)}^{4/3} + \|\chi_{r',r} v\|_{Y(\Rm)}^{4/3}\right) \|\chi_{r',r} (u-v)\|_{Y(\Rm)}. \label{coincidence G 2} 
\end{align}
Here $\chi_{r',r}$ is the characteristic function of the region $\{(x,t): r'+|t|<|x|<r+|t|\}$. We next give an upper bound of $\|\chi_{r',r} (u-v)\|_{Y(\Rm)}$. We follow the same argument as at the beginning of the proof of Lemma \ref{comparison 2}, apply the Strichartz estimates and obtain 
\begin{align*}
 \|\chi_{r',r} (u-v)\|_{Y(\Rm)} \leq \|\chi_{r',r}(u_L-v_L)\|_{Y(\Rm)} + C_2 (\|\chi_{r',r} u\|_{Y(\Rm)}^{4/3} + \|\chi_{r',r} v\|_{Y(\Rm)}^{4/3}) \|\chi_{r',r}(u-v)\|_{Y(\Rm)}. 
\end{align*}
Here $u_L, v_L$ are free waves with initial data $(u_0,u_1)$ and $(v_0,v_1)$. If $r'$ is sufficiently close to $r$, then $\|\chi_{r',r} u\|_{Y(\Rm)}, \|\chi_{r',r} v\|_{Y(\Rm)} \ll 1$. Thus we have 
\begin{equation} \label{coincidence G 3}
 \|\chi_{r',r} (u-v)\|_{Y(\Rm)} \leq 2 \|\chi_{r',r}(u_L-v_L)\|_{Y(\Rm)}.
\end{equation}
By Lemma \ref{u by G} we have
\begin{align*}
 |u_L (x,t)-v_L (x,t)|  =  \left|\frac{1}{|x|^3}\int_{t-|x|}^{t+|x|} (s-t) [G(s)-\tilde{G}(s)] ds\right|.
\end{align*}
If $r'+|t|<|x|<r+|t|$, then $|t\pm |x||>r'$. Thus we have 
\[
 |u_L (x,t)-v_L (x,t)| \leq \frac{1}{|x|^3} \int_{r' < |s|< r} |s-t| |G(s)-\tilde{G}(s)| ds \lesssim_1 \frac{r+|t|}{|x|^3}(r-r')^{1/2}\|G-\tilde{G}\|_{L^2(\{s: |s|>r'\})}.
\]
A straight forward calculation shows ($r' > r/2$)
\[
  \|\chi_{r',r}(u_L-v_L)\|_{Y(\Rm)} \lesssim_1 \left(\frac{r-r'}{r}\right)^{5/7} \|G-\tilde{G}\|_{L^2(\{s: |s|>r'\})}. 
\]
A combination of this inequality with \eqref{coincidence G 1}, \eqref{coincidence G 2} and \eqref{coincidence G 3} shows that if $r'$ is sufficiently close to $r$, then we have 
\[
 \|G-\tilde{G}\|_{L^2(\{s: |s|>r'\})} \lesssim_\gamma \left(\frac{r-r'}{r}\right)^{5/7}  \left(\|\chi_{r',r} u\|_{Y(\Rm)}^{4/3} + \|\chi_{r',r} v\|_{Y(\Rm)}^{4/3}\right) \|G-\tilde{G}\|_{L^2(\{s: |s|>r'\})}. 
\]
Thus if $r'$ is sufficiently close to $r$, then we have 
\[
 \|G-\tilde{G}\|_{L^2(\{s: |s|>r'\})} = 0; \qquad \Rightarrow \qquad G(s)=\tilde{G}(s), \quad \forall |s|>r'.
\]
This immediately gives $u(x,t) = v(x,t)$ as long as $|x|>r'+|t|$. 
\end{proof}

\begin{corollary} \label{asymptotic behaviour with two characteristic numbers}
 Let $\alpha, \beta \in \Rm$ and $u^\alpha$ be the radial weakly non-radiative solution given in Proposition \ref{existence of solution with first alpha}. If $u$ is a radial $R$-weakly non-radiative solution $u$ with first characteristic number $\alpha$ and second characteristic number $\beta$, then $u$ and its initial data $(u_0,u_1)$ satisfies 
 \begin{align*}
  \|\chi_r u\|_{Y(\Rm)} & \lesssim_1 \max\{|\alpha| r^{-1/2}, |\beta|r^{-3/2}\};\\
  \|(u_0,u_1) - (u_0^\alpha, u_1^\alpha) - (\beta |x|^{-3}, 0)\|_{\dot{H}^1 \times L^2 (\{x: |x|>r\})} &\lesssim_1 \max\{c|\beta|^{7/3} r^{-7/2}, c|\alpha|^{4/3}|\beta| r^{-13/6}\}. 
 \end{align*}
 for all $r\geq  \max \{c^{3/2} |\alpha|^2, c^{1/2}|\beta|^{2/3}, R\}$. Here $c = c(\gamma)$ is a constant. 
\end{corollary} 
\begin{proof}
We first show that these inequalities hold for the solution $u$ constructed in Proposition \ref{existence of solution with second beta}. The upper bound of $\|\chi_r u\|$ has been given in \eqref{Y norms of u}. Next we prove the second inequality.  Given $r\geq R=\max \{c^{3/2} |\alpha|^2, c^{1/2}|\beta|^{2/3}\}$, we first choose $G_{0,\beta}$ be radiation profile so that $G_{0,\beta}(s) = 0$ for $|s|>r$ and 
\begin{align*}
 &\int_{-r}^r G_{0,\beta} (s) ds = 0;& &\int_{-r}^r s G_{0,\beta} (s) ds = \beta.& 
\end{align*}
By Lemma \ref{u by G} we have the corresponding initial data $(u_{0,\beta} , u_{1,\beta})$ satisfies 
\begin{equation*} 
 (u_{0,\beta}(x), u_{1,\beta}(x)) = (\beta |x|^{-3}, 0), \qquad |x|>r.
\end{equation*}
The radiation profile of $(u_0,u_1) - (u_0^\alpha, u_1^\alpha) - (u_{0,\beta}, u_{1,\beta})$ is $G-G_{0,\beta}$. Here $G \in X$ is the fixed point we found in the proof of Proposition \ref{existence of solution with second beta}. We then use the identity in the exterior region given above and Corollary \ref{outward decay 2} to conclude that
\begin{align*}
 & \|(u_0,u_1) - (u_0^\alpha, u_1^\alpha) - (\beta |x|^{-3}, 0)\|_{\dot{H}^1 \times L^2 (\{x: |x|>r\})} \\
  & \lesssim_1  r^{-3/2} \left|\int_{-r}^r s(G-G_{0,\beta}) ds\right| + r^{-1/2} \left|\int_{-r}^r (G-G_{0,\beta}) ds\right| + \|G-G_{0,\beta}\|_{L^2(\{s: |s|>r\})}\\
  & \lesssim_1  r^{-3/2} \left|\int_{|s|>r} sG(s) ds\right| + r^{-1/2} \left|\int_{|s|>r} G(s) ds\right| + \|G\|_{L^2(\{s: |s|>r\})}\\
  & \lesssim_1  \max\{c|\beta|^{7/3} r^{-7/2}, c|\alpha|^{4/3}|\beta| r^{-13/6}\}.
\end{align*}
Finally Proposition \ref{uniqueness of solution} implies that two radial weakly non-radiative solutions with the same characteristic numbers must coincide with each other in the overlap part of their exterior regions. Thus the inequalities still hold as long as $r \geq \max \{c^{3/2} |\alpha|^2, c^{1/2}|\beta|^{2/3}, R\}$.
\end{proof}
\section{Dynamics of non-radiative solutions}
Given a radial $R$-weakly non-radiative solution $u(x,t)$ with first characteristic number $\alpha$ and second characteristic number $\beta$, then $v(x,t) \doteq u(x, t+t_0)$ is also a radial $(R+|t_0|)$-weakly non-radiative solution to the time-translated wave equation 
\[
 \partial_t^2 v - \Delta v = F(t+t_0, x, v). 
\]
The non-linear term $F(t+t_0,x,v)$ still satisfies our basic assumption on $F$ with the same norm $\gamma$. A natural question arises that what the first and second characteristic numbers of $v$ are. We first answer the question about the first characteristic number. 

\begin{proposition}
 Let $u$ be a radial weakly non-radiative solution $u(x,t)$ with first characteristic number $\alpha$. Then for all $t_0 \in \Rm$, the radial weakly non-radiative solution $v(x,t) \doteq u(x, t+t_0)$ to the time-translated wave equation given above has the same first characteristic number $\alpha$. 
\end{proposition}
\begin{proof}
 Let us fix $t_0 \in \Rm$. Our main theorem and Proposition \ref{decay of G} give the following estimate if $r > R_1\gg 1$.
 \begin{align*}
  &\|(u_0,u_1) - (0, \alpha |x|^{-3})\|_{\dot{H}^1 \times L^2(\{x: |x|>r\})} \lesssim r^{-7/6}; & &\|\chi_r u\|_{Y(\Rm)} \lesssim r^{-1/2}.&
 \end{align*}
 Then we may apply the Strichartz estimates and finite speed of propagation to obtain 
 \begin{align*}
  \left\|\begin{pmatrix} u(\cdot,t_0)\\ u_t(\cdot,t_0)\end{pmatrix} - \mathbf{S}_L(t_0) \begin{pmatrix} 0\\ \alpha |x|^{-3}\end{pmatrix}\right\|_{\mathcal{H}} & \leq \left\|\begin{pmatrix} u(\cdot,t_0)\\ u_t(\cdot,t_0)\end{pmatrix} - \mathbf{S}_L(t_0) \begin{pmatrix} u_0\\ u_1\end{pmatrix} \right\|_{\mathcal{H}} + \left\|\mathbf{S}_L(t_0) \begin{pmatrix} u_0 - 0\\ u_1 - \alpha |x|^{-3} \end{pmatrix}\right\|_{\mathcal{H}} \\
  & \lesssim_1 \|\chi_r F(t,x,u)\|_{Z(\Rm)} + \left\|\begin{pmatrix} u_0 - 0\\ u_1 - \alpha |x|^{-3} \end{pmatrix}\right\|_{\dot{H}^1 \times L^2(\{x: |x|>r\})}\\
  & \lesssim \|\chi_r u\|_{Y(\Rm)}^{7/3} + r^{-7/6}\\
  & \lesssim r^{-7/6}. 
 \end{align*}
Here $\mathcal{H} = \dot{H}^1 \times L^2 (\{x: |x|>|t_0|+r\})$ and $\mathbf{S}_L(t)$ is the linear wave propagation operator. Next we observe the basic facts 
\[
  \mathbf{S}_L(t_0) \begin{pmatrix} 0\\ \alpha |x|^{-3}\end{pmatrix} = \begin{pmatrix} \alpha t_0 |x|^{-3} \\ \alpha |x|^{-3}\end{pmatrix}.
\]
This gives an inequality 
\begin{equation} \label{estimate of ut L2}
 \|(u(\cdot,t_0),u_t(\cdot,t_0)) - (\alpha t_0 |x|^{-3}, \alpha |x|^{-3})\|_{\dot{H}^1 \times L^2(\{x: |x|>r+|t_0|\})} \lesssim r^{-7/6}, \qquad r > R_1. 
\end{equation}
Please note that the implicit constant in the inequality above does not depend on $t_0$. A basic calculation shows 
\[
 \|\alpha t_0 |x|^{-3}\|_{\dot{H}^1(x:|x|>r+|t_0|)} \lesssim_1 |\alpha||t_0| (r+|t_0|)^{-3/2}
\]
Thus we also have 
\begin{equation*} 
 \|(u(\cdot,t_0),u_t(\cdot,t_0)) - (0, \alpha |x|^{-3})\|_{\dot{H}^1 \times L^2(\{x: |x|>r+|t_0|\})} \lesssim r^{-7/6}, \qquad r > R_1. 
\end{equation*}
This implies that the first characteristic number of $v$ is also $\alpha$. 
\end{proof}

\begin{proposition} \label{dynamics of second}
 Assume that the nonlinear term $F$ is independent of time $t$. Given $\alpha \in \Rm$ and a radial weakly non-radiative solution $u^\alpha$ with first characteristic number $\alpha$, if $u$ is a radial weakly non-radiative solution with first characteristic number $\alpha$ and second characteristic number $\beta$ (with respect to $u^\alpha$), then $v(x,t) \doteq u(x,t+t_0)$ is a radial weakly non-radiative solution to (CP1) with second characteristic number $\beta+\alpha t_0$ (with respect to $u^\alpha$). 
\end{proposition}
\begin{proof}
 Assume that the second characteristic number of $v$ is $\beta'$. Then our main theorem gives that if $r$ is sufficiently large, then
\begin{align*}
 \left\|(u_0,u_1)-(u^\alpha(\cdot,0), u_t^\alpha(\cdot,0)) - (\beta |x|^{-3},0)\right\|_{\dot{H}^1 \times L^2(\{x:|x|>r\})} & \lesssim r^{-13/6};\\
 \left\|(u(\cdot,t_0),u_t(\cdot,t_0))-(u^\alpha(\cdot,0), u_t^\alpha(\cdot,0)) - (\beta' |x|^{-3},0)\right\|_{\dot{H}^1 \times L^2(\{x:|x|>r\})} & \lesssim r^{-13/6}.
\end{align*}
We recall the Hardy's inequality in the exterior region: if $f \in \dot{H}^1 (\Rm^5)$, then
\[
 \|f\|_{L^2(\{x\in \Rm^5: r<|x|<2r\})} \leq 2r\|f(x)/|x|\|_{L^2\{x\in \Rm^5: |x|>r\}} \lesssim_1 r \|f\|_{\dot{H}^1(\{x\in \Rm^5: |x|>r\})}. 
\]
Thus we have 
\begin{align*}
 \left\|u_0 - u^\alpha(\cdot,0) - \beta |x|^{-3}\right\|_{L^2(\{x: |x|>r\})} + \left\|u(\cdot,t_0) - u^\alpha(\cdot,0) - \beta' |x|^{-3}\right\|_{L^2(\{x: |x|>r\})}\lesssim r^{-7/6};
\end{align*}
We may take a difference 
\begin{equation} \label{difference L2}
 \left\|u(\cdot,t_0) - u_0 - (\beta' - \beta) |x|^{-3}\right\|_{L^2(\{x: |x|>r\})} \lesssim r^{-7/6}. 
\end{equation} 
Next we recall the approximation of $u_t$ near the infinity given in \eqref{estimate of ut L2}
\[
 \|u_t(\cdot,t) - \alpha |x|^{-3}\|_{L^2(\{x: |x|>r\})} \lesssim r^{-7/6}, \qquad \forall r\gg 1,\; 0 \leq |t| \leq |t_0|
\]
Thus we may integrate and obtain 
\begin{equation} \label{difference2 L2}
 \|u(\cdot, t_0) - u_0 - \alpha t_0 |x|^{-3}\|_{L^2(\{x: |x|>r\})} \lesssim |t_0| r^{-7/6}, \qquad \forall r \gg 1. 
\end{equation}
We recall the fact $\||x|^{-3}\|_{L^2(\{x:|x|>r\})} \simeq r^{-1/2}$, compare \eqref{difference L2} with \eqref{difference2 L2}, and finally conclude that $\beta' = \beta + \alpha t_0$. 
\end{proof}

\paragraph{Almost periodic solutions} Next we give an application of our characteristic number theory. The soliton-like solutions considered below play an important role in the compactness-rigidity argument, which is a powerful method in the study of global behaviours of solutions to the nonlinear wave or Schr\"{o}dinger equations. Readers may refer to Kenig-Merle \cite{kenig, kenig1} to learn more about the compactness-rigidity argument. This argument consists of two major parts: the compactness part reduces the problem to the study of critical elements, which satisfy strong assumptions on compactness (possibly module symmetries of the equation); The rigidity part deals with these critical elements. Nowadays the compactness has becomes somewhat standard in both wave and Schr\"{o}dinger equations. But the rigidity part does depend on the specific situations. In the radial case we may further reduce the rigidity part into three special situations (see Killip-Tao-Visan \cite{tao}). One of these three situations is the following soliton-like solution. 

\begin{corollary}[Soliton-like solutions] \label{almost periodic solutions}
 Assume that the nonlinear term $F$ is independent of time $t$. Let $u$ be a radial solution to the semi-linear wave equation (CP1) defined for all $t\in \Rm$ so that the orbit 
 \[
  \left\{(u(\cdot,t), u_t(\cdot,t)): t\in \Rm\right\}
 \]
 is a pre-compact subset in $\dot{H}^1 \times L^2(\Rm^5)$, then $u$ is a stationary solution. 
\end{corollary}
\begin{proof}
 It is clear that the translated solution $v^{(t_0)}(x,t) = u(x,t+t_0)$ is a radial non-radiative solution to (CP1) for any given time $t_0 \in \Rm$. Let the first and second characteristic numbers of $u$ be $\alpha$ and $\beta$, respectively. According to Proposition \ref{dynamics of second}, the first and second characteristic numbers of $v^{(t_0)}(x,t)$ are $\alpha$ and $\beta+\alpha t_0$, respectively. We claim that $\alpha = 0$. If this were false, then Corollary \ref{asymptotic behaviour with two characteristic numbers} would give the asymptotic behaviour of $(u(\cdot,t_0), u_t(\cdot, t_0))$. 
 \[
  \|(u(\cdot,t_0), u_t(\cdot, t_0)) - (u_0^\alpha, u_1^\alpha) - ((\beta+\alpha t_0) |x|^{-3}, 0)\|_{\dot{H}^1 \times L^2 (\{x: |x|>r\})} \lesssim_1 c|\beta + \alpha t_0|^{7/3} r^{-7/2}. 
 \]
 Here we assume $t_0$ is sufficiently large and $c^{1/2} |\beta+\alpha t_0|^{2/3} < r < |\beta + \alpha t_0|/|\alpha|$. Here $c = c(\gamma) \gg 1$ is a constant. We choose $r(t_0) = \rho c^{1/2} |\beta+\alpha t_0|^{2/3}$ with a large constant $\rho = \rho(\gamma) > 0$. We have 
 \begin{align*}
  \|((\beta+\alpha t_0) |x|^{-3}, 0)\|_{\dot{H}^1(\{x: |x|>r(t_0)\})} & \gtrsim_1 |\beta+\alpha t_0| r^{-3/2}(t_0) = \rho^{-3/2} c^{-3/4};\\
  c|\beta + \alpha t_0|^{7/3} r^{-7/2} (t_0) & = \rho^{-7/2} c^{-3/4}. 
 \end{align*}
 Therefore we may choose a constant $\rho$ so that $\rho^{-3/2} c^{-3/4} \gg  \rho^{-7/2} c^{-3/4}$. This immediately gives that there exists a small positive constant $\delta \simeq \rho^{-3/2} c^{-3/4}$ so that 
 \[
  \|(u(\cdot,t_0), u_t(\cdot, t_0)) - (u_0^\alpha, u_1^\alpha) \|_{\dot{H}^1 \times L^2 (\{x: |x|>\rho c^{1/2} |\beta+\alpha t_0|^{2/3}\})} \geq \delta.
 \] 
 holds for all sufficiently large $t_0$. This contradicts with the compact property of $(u(\cdot,t_0), u_t(\cdot, t_0))$ because we have 
 \[
  \lim_{t_0\rightarrow \infty} \rho c^{1/2} |\beta+\alpha t_0|^{2/3} = +\infty. 
 \]
 Thus we have $\alpha = 0$. This implies that $v^{(t_0)}(x,t)$ share the same first and second characteristic numbers for all $t_0 \in \Rm$. Thus $v^{(t_0)}(x,t)$ are independent of time $t$. In other words, $u$ is a stationary solution. 
\end{proof}

\noindent This further leads to the following scattering result of defocusing wave equation in dimension 5 with variable coefficient. 

\begin{corollary}
 Consider the defocusing energy-critical wave equation $\partial_t^2 u - \Delta u = - \phi (x) |u|^{4/3} u$. Here the coefficient function $\phi$ is a radial, positive, continuous function so that the limit of $\phi(x)$ as $|x| \rightarrow +\infty$ is well-defined and finite. Then any solution to this equation with radial initial data $(u_0,u_1)\in \dot{H}^1 \times L^2(\Rm^5)$ must be globally defined for all time and scatter in both two time directions. 
\end{corollary}
\begin{proof}
 We follow the compactness-rigidity argument. The compactness part of argument can be found in Section 4 of Shen \cite{compactness35}. This produces a minimal blow-up solution $u$ to the equation above satisfying all of the following conditions if the scattering result above failed. 
\begin{itemize}
 \item The energy $E_\phi (u)$ of $u$ is minimal among all non-scattering solutions; 
 \[
  E_\phi(u) = \int_{\Rm^5} \left(\frac{1}{2} |\nabla u(x,t)|^2 + \frac{1}{2}|u_t(x,t)|^2 + \frac{3}{10}\phi(x) |u(x,t)|^{10/3}\right) dx;
 \]
 \item The solution $u$ is defined for all time $t\in \Rm$ with $\|u\|_{Y([0,+\infty))} = \|u\|_{Y((-\infty,0])} = +\infty$;
 \item The trajectory $\{(u(\cdot, t), u_t(\cdot,t)): t\in \Rm\}$ is pre-compact in the space $\dot{H}^1 \times L^2(\Rm^5)$. 
\end{itemize}
The minimal blow-up solution $u$ here is also a radial function since everything is radial in this situation, although the original compactness argument was given in the non-radial setting. If $\phi(x)$ satisfies an additional assumption 
\[
 \phi(x) - \frac{3}{8} x\cdot \nabla \phi(x) \geq 0, \qquad \forall x\in \Rm^5
\]
then a Morawetz-type estimate does the job. Now we show that the minimal blow-up solution does not exist as long as $\phi(x)$ is positive. According to Corollary \ref{almost periodic solutions}, $u$ must be a stationary solution, i.e. a solution to the elliptic equation 
\[
 -\Delta u = -\phi(x) |u|^{4/3} u. 
\]
Integration by parts then shows 
\[
 \int_{\Rm^5} |\nabla u|^2 dx = - \int_{\Rm^5} \phi(x) |u|^{10/3} dx. 
\]
This gives a contradiction. 
\end{proof} 

\paragraph{Universal profile} We assume that $F$ is independent of time $t$. Given $\alpha,\beta \in \Rm$, we consider the set $\mathcal{R}_{\alpha,\beta}$ of nonnegative real numbers $R$ so that there exists a radial $R$-weakly non-radiative solution to (CP1) with initial data in the energy space, first characteristic number $\alpha$ and second characteristic number $\beta$. We also define $R_{\alpha,\beta} = \inf \mathcal{R}_{\alpha,\beta}$. By the uniqueness of radial weakly non-radiative solutions with given characteristic numbers, we may define a solution $u_{\alpha,\beta}(x,t)$ in the region $\Omega_{R_{\alpha,\beta}}$ so that it is an $R$-weakly non-radiative solution for all $R>R_{\alpha, \beta}$. We call this solution $u_{\alpha,\beta}$ the maximal non-radiative solution to (CP1) with characteristic numbers $\alpha,\beta$. We claim that for any $\alpha \neq 0$ and $\beta_1, \beta_2 \in \Rm$, we have  
\begin{equation} \label{variation of R alpha beta}
  |R_{\alpha,\beta_1} - R_{\alpha, \beta_2}| \leq \frac{|\beta_1-\beta_2|}{|\alpha|}.
\end{equation} 
In fact, $u_{\alpha, \beta_1}(x, t+\frac{\beta_2-\beta_1}{\alpha})$ is a radial weakly non-radiative solution to (CP1) defined in the region $\Omega_{R}$ for any $R > R_{\alpha, \beta_1} + \frac{|\beta_2-\beta_1|}{|\alpha|}$, with characteristic numbers $\alpha$ and $\beta_2$, according to Proposition \ref{dynamics of second}. Thus we have 
\[
 R_{\alpha,\beta_2} \leq R_{\alpha, \beta_1} + \frac{|\beta_2-\beta_1|}{|\alpha|}.
\]
By the symmetric property of $\beta_1, \beta_2$, we obtain \eqref{variation of R alpha beta}. Next we fix $\alpha \in \Rm\setminus \{0\}$ and define a solution $u$ in the following way 
\[
 u(x,t) = u_{\alpha, \alpha t_0} (x,t-t_0), \qquad \forall |x|>R_{\alpha, \alpha t_0}+|t-t_0|, \; t_0 \in \Rm.
\]
This means the time translated version of $u$ defined by $v^{(t_0)}(x,t) = u(x,t+t_0)$ satisfies 
\[
 v^{(t_0)}(x,t) = u_{\alpha, \alpha t_0} (x,t), \qquad \forall |x|>R_{\alpha, \alpha t_0} + |t|;
\]
thus is an $R$-weakly non-radiative solution solution to (CP1) with characteristic numbers $\alpha, \alpha t_0$ for all $R>R_{\alpha, \alpha t_0}$. Namely, all radial weakly non-radiative solutions with first characteristic number $\alpha$ can be viewed as a restriction of $u$, up to a time translation. As a result, we call this solution $u$ the universal profile of radial non-radiative solutions with first characteristic number $\alpha$. Of course we still need to show that $u(x,t)$ is well-defined, since the value of $u$ at $(x,t)$ may be defined for multiple times with different choices of time $t_0$. We fix two times $t_1$, $t_2$ and show that 
\begin{equation} \label{coincidence of solutions}
 u_{\alpha, \alpha t_1} (x, t-t_1) = u_{\alpha, \alpha t_2} (x, t-t_2), \qquad \forall |x|>\max\{R_{\alpha, \alpha t_1}+|t-t_1|, R_{\alpha, \alpha t_2}+|t-t_2|\}.
\end{equation}
\begin{figure}[h]
 \centering
 \includegraphics[scale=0.8]{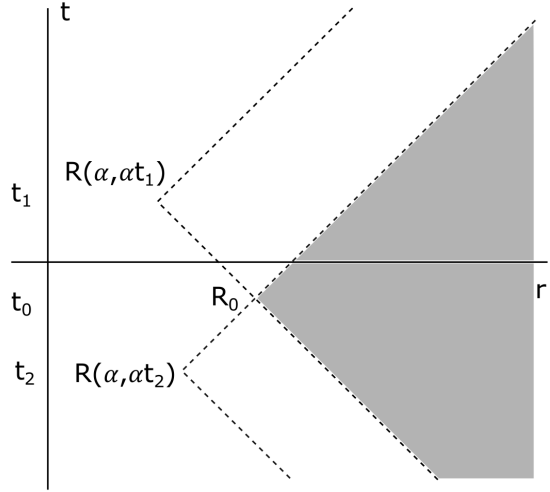}
 \caption{Overlap region} \label{overlap region}
\end{figure}
We may rewrite the overlap region, which contains all $(x,t)$ satisfying the inequality above, in the form of 
\[
 \{(x,t)\in \Rm^5 \times \Rm: |x|>R_0+|t-t_0|)\},
\]
as illustrated in figure \ref{overlap region}. The inequality \eqref{variation of R alpha beta} implies that we always have 
\begin{align*}
 &R_0 = R_{\alpha, \alpha t_1} + |t_1 -t_0|;& &R_0 = R_{\alpha, \alpha t_2} + |t_2 -t_0|.&
\end{align*}
We apply Proposition \ref{dynamics of second} and obtain that $u_{\alpha, \alpha t_i}(x, t + t_0-t_i)$ is a radial weakly non-radiative solution to (CP1) in the exterior region $\Omega_R$ for any $R > R_{\alpha, \alpha t_i} + |t_i-t_0| = R_0$, with characteristic numbers $\alpha$ and $t_0 \alpha$. By the uniqueness of radial weakly non-radiative solutions with given characteristic numbers, we have 
\[
 u_{\alpha, \alpha t_1}(x, t + t_0-t_1) = u_{\alpha, \alpha t_2}(x, t + t_0-t_2), \qquad \forall |x|>R_0 + |t|.
\]
A substitution of $t$ by $t-t_0$ then gives the identity \eqref{coincidence of solutions}. In summary we have 

\begin{proposition} \label{universal profile alpha}
 Assume that $F$ is independent of time $t$. Given $\alpha \in \Rm \setminus \{0\}$, there exists a radial function $u^\alpha (x,t)$ defined in the region $\{(x,t)\in \Rm^5 \times \Rm: |x|>R_{\alpha, \alpha t}\}$ so that 
 \begin{itemize} 
  \item The radius $R_{\alpha, \alpha t}$ satisfies $|R_{\alpha, \alpha t_1}- R_{\alpha, \alpha t_2}| \leq |t_1-t_2|$ for all $t_1,t_2\in \Rm$.
  \item If $t_0\in \Rm$ and $R > R_{\alpha, \alpha t_0}$, then $u^\alpha (x,t+t_0)$ is a radial $R$-weakly non-radiative solution to (CP1) with characteristic numbers $\alpha, \alpha t_0$. 
  \item Conversely, if $v$ is a radial $R$-weakly non-radiative solution to (CP1) with characteristic numbers $\alpha, \alpha t_0$, then $R \geq R_{\alpha, \alpha t_0}$ and $v(x,t) = u^\alpha (x,t_0+t)$ holds as long as $|x|>R+|t|$. 
 \end{itemize}
\end{proposition}

\begin{remark}
 If we consider the focusing/defocusing wave equation $\partial_t^2 u - \Delta u = \pm |u|^{4/3} u$, then any radial weakly non-radiative solution $v$ to (CP1) with a nonzero first characteristic number can be viewed as a restriction of $u^1$ defined above, up to a time translation, a sign and/or a dilation.  Namely there exists $(\lambda, t_0, \zeta) \in \Rm^+ \times \Rm \times \{-1,1\}$ so that the following identity holds in the corresponding exterior region
 \[
  v(x,t) = \zeta \frac{1}{\lambda^{3/2}} u^1 \left(\frac{x}{\lambda}, \frac{t+t_0}{\lambda} \right). 
 \]
\end{remark}
\section*{Acknowledgement}
The second author is financially supported by National Natural Science Foundation of China Project 12071339.

\end{document}